\documentclass[final]{elsarticle}
\usepackage[english]{babel}
\usepackage{lineno,hyperref}

\usepackage[T1]{fontenc}
\usepackage[latin9]{inputenc}
\usepackage{geometry}
\usepackage{color}
\usepackage{float}
\usepackage{amsmath}
\usepackage{amssymb}
\usepackage{graphicx}
\usepackage{esint}

\makeatletter

\providecommand{\tabularnewline}{\\}

\@ifundefined{showcaptionsetup}{}{%
 \PassOptionsToPackage{caption=false}{subfig}}
\usepackage{subfig}
\makeatother

\usepackage{babel}

\modulolinenumbers[5]

\journal{Journal of XXX}

\begin{document}

\begin{frontmatter}
\date{April, 25, 2018}
\title{A Numerical Study of Steklov Eigenvalue Problem via Conformal Mapping\tnoteref{mytitlenote}}
\tnotetext[mytitlenote]{The work of Chiu-Yen Kao is partially supported
by a Collaboration Grant for Mathematicians 514210 from the Simons
Foundation.}


\author[mymainaddress]{Weaam Alhejaili}

\author[mysecondaryaddress]{Chiu-Yen Kao\corref{mycorrespondingauthor}}
\cortext[mycorrespondingauthor]{Corresponding author}
\ead{ckao@cmc.edu}

\address[mymainaddress]{E-mail: weaam.alhejaili@cgu.edu; Address: Institute of Mathematical
Sciences, Claremont Graduate University, Claremont, CA 91711}
\address[mysecondaryaddress]{E-mail: ckao@cmc.edu; Address: Department of Mathematical Sciences,
Claremont McKenna College, Claremont, CA 91711}

\begin{abstract}
In this paper, a spectral method based on conformal mappings is proposed
to solve Steklov eigenvalue problems and their related shape optimization
problems in two dimensions. To apply spectral methods, we first reformulate
the Steklov eigenvalue problem in the complex domain via conformal
mappings. The eigenfunctions are expanded in Fourier series so the
discretization leads to an eigenvalue problem for coefficients of
Fourier series. For shape optimization problem, we use the gradient
ascent approach to find the optimal domain which maximizes $k-$th
Steklov eigenvalue with a fixed area for a given $k$. The coefficients
of Fourier series of mapping functions from a unit circle to optimal
domains are obtained for several different $k$.\end{abstract}

\begin{keyword}
\texttt{Steklov eigenvalues, extremal eigenvalue problem,
shape optimization, spectral method, conformal mapping}
\MSC[2010] 35P15 \sep 49Q10 \sep 65N25 \sep 65N35
\end{keyword}

\end{frontmatter}

\linenumbers

\section{Introduction }

The second order Steklov eigenvalue problem satisfies

\begin{equation}
\left\{ \begin{array}{ccc}
\triangle u(\mathbf{x})=0 & \text{in} & \Omega,\\
\partial_{n}u=\lambda u & \text{on} & \partial\Omega,
\end{array}\right.\label{eq:Steklov eigenvalue}
\end{equation}
where $\triangle$ is the Laplace operator acting on the function
$u(\mathbf{x})$ defined on $\Omega\subset\mathbb{R}^{N}$, $\lambda$
is the corresponding eigenvalue, and $\partial_{n}$ is the outward
normal derivative along the boundary $\partial\Omega.$ This problem
is a simplified version of the mixed Steklov problem which was used
to obtain the sloshing modes and frequencies. The spectral geometry
of the Steklov problem has been studied for a long time. See a recent
review article on American Mathematical Society (AMS) notice \cite{kuznetsov2014legacy}
and the references therein.\textcolor{black}{{} In 2012, Krechetnikov
and Mayer were awarded the Ig Noble prize for fluid dynamics for their
work on the dynamic of liquid sloshing. In \cite{mayer2012walking},
they studied the conditions under which coffee spills for various
walking speeds based on sloshing modes \cite{ibrahim2005liquid}. }

The Steklov problem (\ref{eq:Steklov eigenvalue}) has a countable
infinite set of eigenvalues which are greater than or equal to zero.
We arrange them as $0=\lambda_{0}(\Omega)<\lambda_{1}(\Omega)\leq\lambda_{2}(\Omega)\leq\cdots\leq\lambda_{k}(\Omega)\leq\cdots\rightarrow\infty$
and denote $u_{k}\in H^{1}(\Omega)$ as the corresponding eigenfunction.
The Weyl\textquoteright s law for Steklov eigenvalues states that
$\lambda_{k}\sim2\pi\left(\frac{k}{\mid\mathbb{B}^{N-1}\mid\mid\partial\Omega\mid}\right)^{\frac{1}{N-1}}$
where $\mathbb{B}^{N-1}$ is the unit ball in $\mathbb{R}^{N-1}$.
The variational characterization of the eigenvalues is given by 
\begin{equation}
\lambda_{k}(\Omega)=\underset{v\in H^{1}}{\min}\left\{ \frac{\int_{\Omega}\left|\nabla v\right|^{2}dx}{\int_{\partial\Omega}v^{2}ds}:\int_{\partial\Omega}vu_{i}=0,\;i=0,\ldots,k-1\right\} .\label{eq: variational formulation of lambda}
\end{equation}

In 1954, Weinstock proved that the disk maximizes the first non-trivial
Steklov eigenvalue $\lambda_{1}$ among simply-connected planar domains
with a fixed perimeter \cite{weinstock1954inequalities,girouard2010shape}.
Furthermore, the $k$-th eigenvalue $\lambda_{k}$ for a simply-connected
domain with a fixed perimeter is maximized in the limit by a sequence
of simply-connected domains degenerating to the disjoint union of
\textcolor{black}{$k$ identical disks for any $k\geq1$ \cite{girouard2010hersch}.}
It remains an open question for non-simply-connected bounded planar
domains \cite{girouard2017spectral}.\textcolor{black}{{} Furthermore,
the existence of the optimal shapes that maximized the Steklov eigenvalues
was proved in \cite{bogosel2017steklov} recently. }

Several different numerical approaches were proposed to solve Steklov
eigenvalue problem \cite{bogosel2016method,akhmetgaliyev2017computational}
and Wentzell eigenvalue problem \cite{bogosel2016method} which has
slightly different boundary conditions. The methods of fundamental
solutions were used in \cite{bogosel2016method} to compute Steklov
spectrum and a theoretical error bound were derived. In \cite{akhmetgaliyev2017computational},
the authors used a boundary integral method with a single layer potential
representation of eigenfunction. Both methods can possibly achieve
spectral convergence. Furthermore, they both studied maximization
of $\lambda_{k}$ among star-shaped domains with a fixed area \cite{akhmetgaliyev2017computational,bogosel2016method}.

Mixed boundary problems were solved in \cite{andreev2004isoparametric}
and \cite{mora2015virtual} via isoparametric finite element method
and the virtual element method, respectively. The error estimates
for eigenvalues and eigenfunctions were derived. Another type of Steklov
problem which is formulated as
\[
\left\{ \begin{array}{ccc}
-\triangle u(\mathbf{x})+u(\mathbf{x})=0 & \text{in} & \Omega,\\
\partial_{n}u=\lambda u & \text{on} & \partial\Omega,
\end{array}\right.
\]
was studied numerically in \cite{bi2016adaptive,xie2013type,bi2011two,li2011two}.
In \cite{bonder2007optimization}, the authors look for a subset $A\subset\Omega$
that minimizes the first Steklov-like problem 
\[
\left\{ \begin{array}{ccc}
-\triangle u(\mathbf{x})+u(\mathbf{x})=0 & \text{in} & \Omega\backslash\bar{A},\\
\partial_{n}u=\lambda u & \text{on} & \partial\Omega,\\
u=0 & \text{on} & \partial A,
\end{array}\right.
\]
by using\textcolor{black}{{} an algorithm based on finite element methods
and shape derivatives.}\textcolor{magenta}{{} }Furthermore, finite element
methods have been also applied to the nonlinear Steklov eigenvalue
problems \cite{kumar2010simulation} and methods of fundamental solutions
were proposed lately to find a convex shape that has the least biharmonic
Steklov eigenvalue \cite{antunes2013convex}.

The aim of this paper is two-fold. First, we develop numerical approaches
to solve the forward problem of Steklov eigenvalue problem by using
spectral methods for complex formulations via conformal mapping approaches
\cite{brown1996complex,kao2010moving} for any given simply-connected
planar domain. Second, we aim to find the maximum value of $\lambda_{k}$
with a fixed area among simply-connected domains via the gradient
ascent approach. To find optimal domains, we start with a chosen initial
domain of any shape and deform the domain with the velocity which
is obtained by calculating the shape derivative of $\lambda_{k}\sqrt{|\Omega|}$
and choose the ascent direction. In the complex formulation, the deforming
domain is mapped to a fixed unit circle which allows spectral methods
to solve the problem efficiently. 

In Section 2, we briefly review the derivation of Steklov eigenvalue
problem. The formulations of Steklov eigenvalue problem in $\mathbb{R}^{2}$
and $\mathbb{C}$ are described in Sections 3 and 4, respectively.
Some known analytical solutions are provided and optimization of $k-$th
Steklov eigenvalue $\lambda_{k}$ is formulated. In Section 5, computational
methods are described and numerical experiments are presented. The
summary and discussion are given in Section 6. 

\section{The derivation of Steklov problem}

Let us briefly review the derivation of Steklov eigenvalue problem
coming from the sloshing model which neglects the surface tension
\cite{ibrahim2005liquid}. Consider the sloshing problem in a three-dimensional
simply-connected container filled with inviscid, irrotational, and
incompressible fluid. Choose Cartesian coordinates $(x,y,z)$ so that
the mean free surface lies in the $(x,y)$-plane and the $z$-axis
is directed upwards. Denote $\tilde{F}$ as the free fluid surface
and $B$ as the rigid bottom of the container. The governing equations
in $\tilde{\Omega}$ of the sloshing model are
\[
\begin{array}{rrcl}
\text{Navier-Stokes equation:} & \frac{\partial\mathbf{V}}{\partial t}+(\mathbf{\boldsymbol{V}}\cdot\nabla)\mathbf{V} & = & -\frac{1}{\rho}\nabla p-\nabla(gz)\\
\text{irrotational flow:} & \nabla\times\mathbf{V} & = & 0\\
\text{incompressible fluid:} & \nabla\cdot\boldsymbol{V} & = & 0\\
\text{velocity potential:} & \boldsymbol{V} & = & \nabla\tilde{\Phi}
\end{array}
\]
where $\boldsymbol{V}(x,y,z,t)$ is the fluid velocity, $\rho$ is
the density, $p$ is the pressure, $g$ is the gravity, and $\tilde{\Phi}(x,y,z,t)$
is the velocity potential. The last two equations lead to Laplace's
equation

\[
\begin{array}{ccc}
\triangle\tilde{\Phi}=0 & \text{in} & \tilde{\Omega}.\end{array}
\]
The no penetration boundary condition at the rigid bottom of the container
is 

\begin{equation}
\begin{array}{ccc}
\nabla\tilde{\Phi}\cdot\hat{n}_{B}=0 & \text{on} & B\end{array}\label{eq:bounary condtion on B}
\end{equation}
where $\hat{n}_{B}$ is the outward unit normal to the boundary $B$
and the dynamic boundary condition at the free surface $z=\tilde{\gamma}(x,y,t)$
is
\begin{equation}
\tilde{\gamma}_{t}+\nabla\tilde{\Phi}\cdot\nabla(\tilde{\gamma}-z)=0.\label{eq: boundary condition on F}
\end{equation}
Rewriting the Navier-Stokes equation in terms of $\tilde{\Phi}$ and
using 
\[
(\mathbf{\boldsymbol{V}}\cdot\nabla)\mathbf{V}=\frac{1}{2}\nabla\left|\mathbf{\boldsymbol{V}}\right|^{2}-\mathbf{V}\times(\nabla\times\mathbf{V})=\frac{1}{2}\nabla\left|\mathbf{\boldsymbol{V}}\right|^{2},
\]
we obtain the Bernoulli's equation

\begin{equation}
\nabla\left(\tilde{\Phi}_{t}+\frac{p}{\rho}+\frac{1}{2}\left|\mathbf{\nabla}\tilde{\Phi}\right|^{2}+gz\right)=0.\label{eq:bernouill's eq}
\end{equation}
Thus

\begin{equation}
\tilde{\Phi}_{t}+\frac{p}{\rho}+\frac{1}{2}\left|\mathbf{\nabla}\tilde{\Phi}\right|^{2}+gz=A(t)\label{eq: integarting Bernulii's}
\end{equation}
where $A(t)$ is an arbitrary function of $t$. By using the condition
that the pressure $p$ at the free surface equals to the ambient pressure
$p_{atm}$ and choosing $A(t)=\frac{p_{atm}}{\rho},$ we then have

\[
\tilde{\Phi}_{t}+\frac{1}{2}\left|\mathbf{\nabla}\tilde{\Phi}\right|^{2}+gz=0.
\]
Therefore, we obtain the following partial differential equations
\begin{equation}
\begin{array}{rclccc}
\triangle\tilde{\Phi} & = & 0 &  & \text{in} & \tilde{\Omega},\\
\nabla\tilde{\Phi}\cdot\hat{n_{B}} & = & 0 &  & \text{on} & B,\\
\tilde{\gamma}_{t}+\nabla\tilde{\Phi}\cdot\nabla(\tilde{\gamma}-z) & = & 0 &  & \text{on} & \tilde{F},\\
\tilde{\Phi}_{t}+\frac{1}{2}\left|\nabla\tilde{\Phi}\right|^{2}+gz & = & 0 &  & \text{on} & \tilde{F}.
\end{array}\label{eq:nonlinear_partial_diffrential}
\end{equation}

Assuming the liquid motion is of small amplitude $z=\tilde{\gamma}(x,y,t)$
from the undisturbed free surface $z=0$, we consider the following
asymptotic expansion:
\[
\begin{array}{rcl}
\tilde{\Phi}(x,y,z,t) & = & \Phi_{0}+\epsilon\hat{\Phi}(x,y,z,t),\\
\tilde{\gamma}(x,y,t) & = & \gamma_{0}+\epsilon\hat{\gamma}(x,y,t),
\end{array}
\]
where $\Phi_{0}$ is a constant velocity potential, $\gamma_{0}=0$,
$\hat{\Phi}(x,y,z,t)$ and $\hat{\gamma}(x,y,t)$ represent perturbations,
and $\epsilon>0$ is a small parameter. Substituting these expansions
in (\ref{eq:nonlinear_partial_diffrential}) gives
\begin{equation}
\begin{array}{rclccc}
\triangle\hat{\Phi} & = & 0 &  & \text{in} & \tilde{\Omega},\\
\nabla\hat{\Phi}\cdot\hat{n_{B}} & = & 0 &  & \text{on} & B,\\
\hat{\gamma}_{t}+\nabla\hat{\Phi}\cdot\nabla(\epsilon\hat{\gamma}-z) & = & 0 &  & \text{on} & \tilde{F},\\
\hat{\Phi}_{t}+\epsilon\frac{1}{2}\left|\mathbf{\nabla}\hat{\Phi}\right|^{2}+g\hat{\gamma} & = & 0 &  & \text{on} & \tilde{F}.
\end{array}\label{eq:undestirb_equations}
\end{equation}
It is well known that the time harmonic solutions of (\ref{eq:undestirb_equations})
with angular frequency $\alpha$ and phase shift $\sigma$ are given
by 
\[
\hat{\Phi}(x,y,z,t)=U(x,y,z)\text{cos\ensuremath{(\alpha t+\sigma),}}
\]
\[
\hat{\gamma}(x,y,t)=\mu(x,y)\text{sin\ensuremath{(\alpha t+\sigma).}}
\]
where $U(x,y,z)$ is the sloshing velocity potential and $\mu(x,y)$
is the sloshing height. Substitute these expansions into (\ref{eq:undestirb_equations}),
transform the boundary conditions on $\tilde{F}$ to $F$ and the
domain $\tilde{\Omega}$ to $\Omega$ by using Taylor expansion about
$z=0$, and ignore high order terms. We then obtain
\[
\begin{array}{rclccc}
\triangle U & = & 0 &  & \text{in} & \Omega,\\
\nabla U\cdot\hat{n_{B}} & = & 0 &  & \text{on} & B,\\
U_{Z} & = & \alpha\mu &  & \text{on} & F,\\
\mu & = & \alpha\frac{U}{g} &  & \text{on} & F.
\end{array}
\]
Thus, we obtain the mixed Steklov eigenvalue problem 
\[
\begin{array}{rclccc}
\triangle U & = & 0 &  & \text{in} & \Omega,\\
\nabla U\cdot\hat{n_{B}} & = & 0 &  & \text{on} & B,\\
U_{Z} & = & \lambda U &  & \text{on} & F,
\end{array}
\]
where $\lambda=\alpha^{2}/g.$

When $B$ is an empty set, the mixed Steklov eigenvalue problem is
reduced to the classical Steklov eigenvalue problem (\ref{eq:Steklov eigenvalue}).
The Steklov spectrum satisfying (\ref{eq:Steklov eigenvalue}) is
also of fundamental interest as it coincides with the spectrum of
the Dirichlet-to-Neumann operator $\Gamma:H^{\frac{1}{2}}(\partial\Omega)\rightarrow H^{-\frac{1}{2}}(\partial\Omega)$,
given by the formula $\Gamma u=\partial_{n}(\mathbb{H}u)$, where
$\mathbb{H}u$ denotes the unique harmonic extension of $u\in H^{\frac{1}{2}}(\partial\Omega)$
to $\Omega.$

\section{Steklov Eigenvalue Problems on $\Omega\subset\mathbb{R}^{2}$}

In this section, we discuss some known analytical solutions of Steklov
eigenvalue problems on simple geometric shapes and formulate the maximization
of Steklov eigenvalue with a fixed area constraint. 

\subsection{Some Known Analytical Solutions}

\subsubsection{On a Circular Domain }

\begin{figure}[H]
\begin{centering}
\includegraphics[scale=0.25]{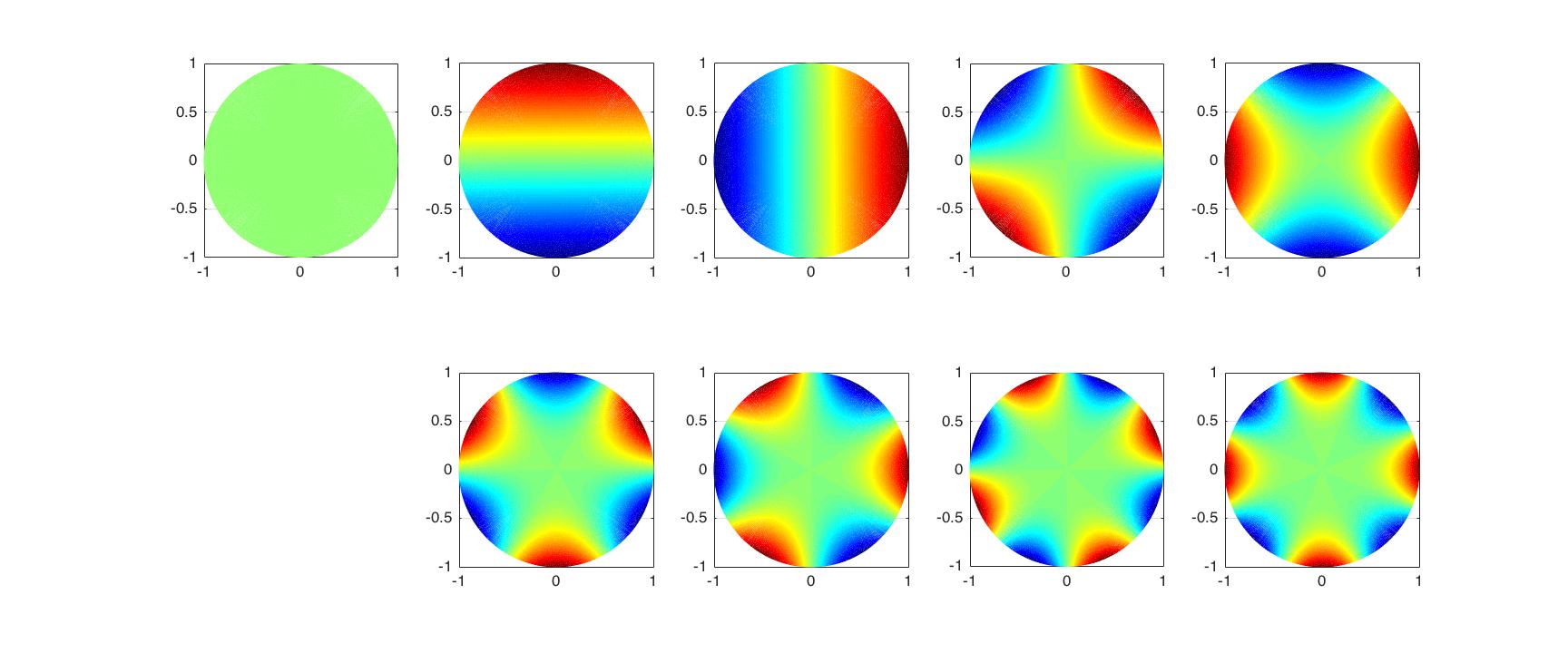}
\par\end{centering}
\caption{The first nine Steklov eigenfunctions on a unit circle.\label{fig:steklov-eigenfunction-for}}

\end{figure}

By using the method of separation of variables, it is well known that
the Steklov eigenvalues of a unit circle $\Omega$ are given by 
\[
0,1,1,2,2,\cdots,k,k,\cdots
\]
where $\lambda_{2k}=\lambda_{2k-1}=k$ has multiplicity 2 and their
corresponding eigenfunctions are
\[
u_{2k}=r^{k}\cos(k\theta),\;u_{2k-1}=r^{k}\sin(k\theta).
\]
The first nine eigenfunctions are shown in Figure \ref{fig:steklov-eigenfunction-for}.

\subsection{\label{subsec:On-an-annulus}On an Annulus}

When $\Omega=B(0,1)\setminus B(0,\epsilon)$, the Steklov eigenvalues
can be found via the method of separation of variables \cite{girouard2017spectral}.
The only eigenfunction which is radial independent satisfies 
\[
u(r)=(\frac{-(1+\epsilon)}{\epsilon\ln\epsilon})\ln(r)+1,
\]
and the corresponding eigenvalue is 
\[
\lambda=\frac{1+\epsilon}{\epsilon}\ln(1/\varepsilon).
\]

The rest of the eigenfunctions are of the form 
\begin{equation}
u_{k}(r,\theta)=(Ar^{k}+Br^{-k})H(k\theta),\quad k\in\mathbb{N}\label{eq:eigenfunction}
\end{equation}
where $A$ and $B$ are constants and $H(k\theta)=\cos(k\theta)$
or $H(k\theta)=\sin(k\theta)$. The boundary conditions become

\begin{eqnarray}
\frac{\partial}{\partial_{r}}u_{k}(1,\theta) & = & \lambda\,u_{k}(1,\theta),\nonumber \\
\frac{\partial}{\partial_{r}}u_{k}(\epsilon,\theta) & = & -\lambda\,u_{k}(\epsilon,\theta),\label{eq:ellipse eigenvalue equation}
\end{eqnarray}
which can be simplified to the following system

\[
\left[\begin{array}{cc}
\lambda\epsilon^{k}+k\epsilon^{k-1} & \lambda\epsilon^{-k}-k\epsilon^{-k-1}\\
\lambda-k & \lambda+k
\end{array}\right]\left[\begin{array}{c}
A\\
B
\end{array}\right]=\left[\begin{array}{c}
0\\
0
\end{array}\right].
\]
To obtain nontrivial solutions, the determinant of the matrix needs
to be zero. Thus Steklov eigenvalues are determined by the roots of
the following polynomial 

\begin{equation}
p_{k}(\lambda)=\lambda^{2}-\lambda k\left(\frac{\epsilon+1}{\epsilon}\right)\left(\frac{1+\epsilon^{2k}}{1-\epsilon^{2k}}\right)+\frac{1}{\epsilon}k^{2},\quad k\in\mathbb{N}.\label{eq:poly eqn}
\end{equation}
Note that every root corresponds to a double eigenvalue. If $\epsilon>0$
is smaller enough, for $k=1$, we get the smallest eigenvalue 
\[
\lambda_{1}(\Omega)=\frac{1}{2\epsilon}\frac{1+\epsilon^{2}}{1-\epsilon}\left(1-\sqrt{1-4\epsilon\left(\frac{1-\epsilon}{1+\epsilon^{2}}\right)^{2}}\right).
\]

\subsection{Shape Optimization }

It follows from (\ref{eq: variational formulation of lambda}) that
the Steklov eigenvalues satisfy the homothety property $\lambda_{k}(t\Omega)=t^{-1}\lambda_{k}(\Omega).$
Instead of fixing the perimeter or the area, one can consider the
following shape optimization problems 
\begin{equation}
\lambda_{k}^{L\star}=\underset{\Omega\subset\mathbb{R}^{2}}{\text{max}}\lambda_{k}^{L}(\Omega)\qquad\text{where \ensuremath{\lambda_{k}^{L}(\Omega)}=\ensuremath{\lambda_{k}(\Omega)\left|\partial\Omega\right|}}\label{eq:optmization length-1}
\end{equation}
and
\begin{equation}
\lambda_{k}^{A\star}=\underset{\Omega\subset\mathbb{R}^{2}}{\text{max}}\lambda_{k}^{A}(\Omega)\qquad\text{where \ensuremath{\lambda_{k}^{A}(\Omega)}=\ensuremath{\lambda_{k}(\Omega)\sqrt{\left|\Omega\right|}}}.\label{eq:optmization area-1}
\end{equation}
\textcolor{black}{As mentioned in the Introduction section, the perimeter
eigenvalue problem (\ref{eq:optmization length-1}) is known analytically
for simply-connected domains. Thus, we focus only on normalized eigenvalue
with respect to the area as described in (\ref{eq:optmization area-1}). }

\subsubsection{On an Annulus}

\begin{figure}[H]
\begin{centering}
\includegraphics[scale=0.6]{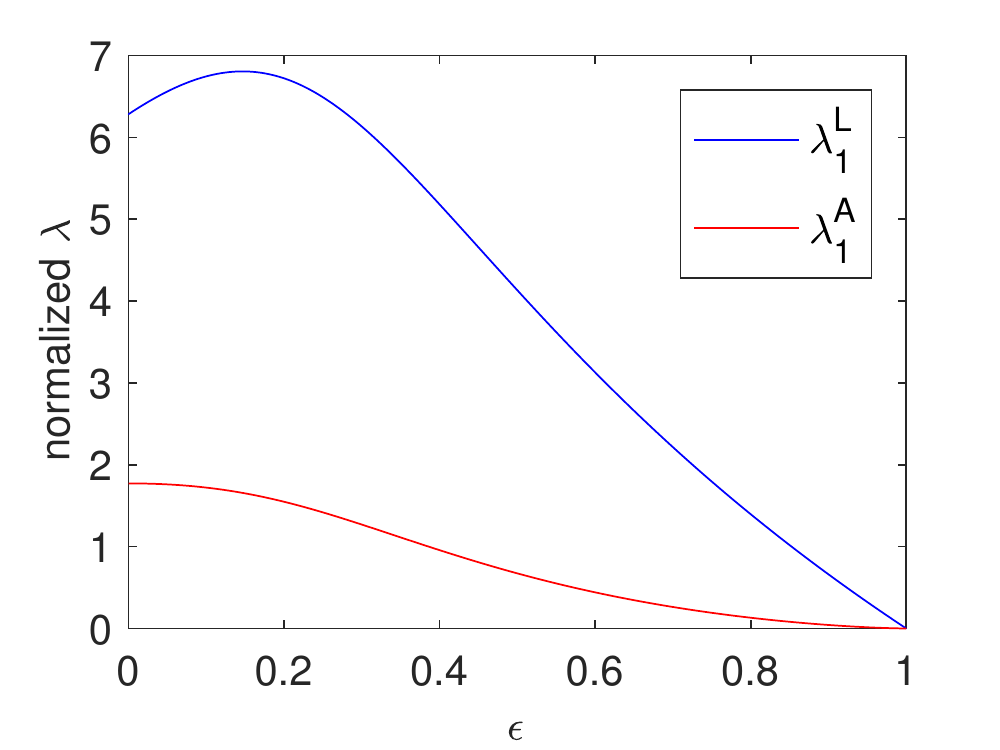}
\par\end{centering}
\begin{centering}
\caption{The perimeter- and area-normalized eigenvalue, $\lambda_{1}^{L}$
and $\lambda_{1}^{A}$, on an annulus, respectively. \label{fig:the-normalized-eigenvalue}}
\par\end{centering}
\end{figure}

In Section \ref{subsec:On-an-annulus} we get $\lambda_{1}(\Omega)$
on an annulus $\Omega=B(0,1)\setminus B(0,\epsilon)$. Thus, $\lambda_{1}^{L}=\lambda_{1}[2\pi(1+\epsilon)]$
is the normalized first eigenvalue with respect to the perimeter of
the domain $\Omega$. The perimeter normalized eigenvalue is not a
monotone function in $\epsilon$ and it reac\textcolor{black}{hes
the maximum value $6.8064$ when $\epsilon=\epsilon^{*}\approx0.1467$}
\cite{girouard2017spectral}\textcolor{black}{{} as shown in Figure
\ref{fig:the-normalized-eigenvalue}. On the other hand $\lambda_{1}^{A}=\lambda_{1}[\sqrt{\pi(1-\epsilon^{2})}]$
is the normalized first eigenvalue with respect to the area of the
domain $\Omega$ which turns out to be a monotone decreasing function
}in $\epsilon$ and it reac\textcolor{black}{hes the maximum value
$\sqrt{\pi}$ when $\epsilon=0$ as shown in Figure \ref{fig:the-normalized-eigenvalue}.}

\subsubsection{Shape derivative}

Here we review the concept of the shape derivative. For more details,
we refer the readers to  \cite{sokolowski1992introduction}.

\textbf{Definition:} Let $\Omega\subset\mathbb{R}^{N}$ and $J$ be
a functional on $\Omega\mapsto J(\Omega)$. Consider the perturbation
$x\in\Omega\rightarrow x+tV\in\Omega_{t}$ where $V$ is a vector
field. Then the shape derivative of the functional $J$ at $\Omega$
in the direction of a vector field $V$ is given by 

\begin{equation}
dJ(\Omega;V)=\underset{t\downarrow0}{lim}\frac{J(\Omega_{t})-J(\Omega)}{t}.\label{eq:shape drivative-1}
\end{equation}

\noindent In \cite{akhmetgaliyev2017computational}, the shape derivative
of Steklov eigenvalue is given by the following proposition.

\textbf{Proposition:} Consider the perturbation $x\mapsto x+tV$ and
denote $c=V\cdot\hat{n}$ where $\hat{n}$ is the outward unit normal
vector. Then a simple (unit-normalized) Steklov eigenpair $\left(\lambda,u\right)$
satisfies the perturbation formula
\begin{equation}
\lambda(\Omega)^{\prime}=\int_{\partial\Omega}(\left|\nabla u\right|^{2}-2\lambda^{2}u^{2}-\lambda\kappa u^{2})c\,ds\label{eq:shape drivative of eigenvalue}
\end{equation}
where $\kappa$ is the mean curvature.

Proof. By using the variational formulation (\ref{eq: variational formulation of lambda})
of eigenvalue and normalizing the eigenfunction by
\begin{equation}
\int_{\partial\varOmega}u^{2}ds=1,\label{eq:normalized equation-1}
\end{equation}
 we have 
\[
\lambda(\Omega)=\int_{\Omega}\left|\nabla u\right|dx.
\]
Now denote the shape derivative by the prime, thus 
\[
\begin{array}{rclcc}
\lambda^{\prime}(\Omega) & = & \left(\int_{\Omega}\left|\nabla u\right|^{2}dx\right)^{\prime} &  & (\text{shape derivative)}\\
\\
 & = & \int_{\Omega}(\left|\nabla u\right|^{2})^{\prime}dx+\int_{\partial\Omega}\left|\nabla u\right|^{2}V\cdot n\,ds\\
\\
 & = & \int_{\Omega}(\nabla u\cdot\nabla u)^{\prime}dx+\int_{\partial\Omega}\left|\nabla u\right|^{2}V\cdot n\,ds\\
\\
 & = & 2\int_{\Omega}\nabla u\cdot(\nabla u)^{\prime}dx+\int_{\partial\Omega}\left|\nabla u\right|^{2}c\,ds\\
\\
 & = & -2\int_{\Omega}(\triangle u)u^{\prime}dx+2\int_{\partial\Omega}u_{n}u^{\prime}ds+\int_{\partial\Omega}\left|\nabla u\right|^{2}c\,ds &  & (\text{Green's identity)}\\
\\
 & = & 2\lambda\int_{\partial\Omega}uu{}^{\prime}ds+\int_{\partial\Omega}\left|\nabla u\right|^{2}c\,ds &  & \text{( Equation (\ref{eq:Steklov eigenvalue}) )}
\end{array}
\]
Now applying the shape derivative to (\ref{eq:normalized equation-1}),
we get
\[
\int_{\partial\Omega}uu^{\prime}\,ds=-\int_{\partial\Omega}\left(uu_{n}+\frac{\kappa}{2}\right)u^{2}c\,ds=-\int_{\partial\Omega}\left(\lambda+\frac{\kappa}{2}\right)u^{2}c\,ds
\]
Therefore, we get (\ref{eq:shape drivative of eigenvalue}) where
$\kappa$ is the mean curvature. 

Now consider the optimization problem (\ref{eq:optmization area-1})
and use the shape derivative of $\lambda$, we get 
\begin{equation}
\begin{array}{rl}
\left(\lambda_{k}^{A}(\Omega)\right)^{\prime} & =\left(\lambda_{k}(\Omega)\cdot\sqrt{\left|\Omega\right|}\right)^{\prime}\\
\\
 & =\lambda_{k}^{\prime}(\Omega)\sqrt{\left|\Omega\right|}+\lambda_{k}(\Omega)\frac{1}{2\sqrt{\left|\Omega\right|}}\Omega^{\prime}\\
\\
 & =\sqrt{\left|\Omega\right|}\int_{\partial\Omega}\left(\left|\nabla u\right|^{2}-2\lambda^{2}u^{2}-\lambda\kappa u^{2}\right)c\,ds+\lambda_{k}(\Omega)\frac{1}{2\left|\Omega\right|}\int_{\partial\Omega}c\,ds\\
\\
 & =\sqrt{\left|\Omega\right|}\int_{\partial\Omega}\left(\left(\left|\nabla u\right|^{2}-2\lambda^{2}u^{2}-\lambda\kappa u^{2}\right)+\lambda_{k}(\Omega)\frac{1}{2\left|\Omega\right|}\right)c\,ds.
\end{array}\label{eq}
\end{equation}
Thus the normalized velocity for the ascent direction can be chosen
as
\begin{equation}
c=V_{n}=\left(\left|\nabla u\right|^{2}-2\lambda^{2}u^{2}-\lambda\kappa u^{2}\right)+\lambda_{k}(\Omega)\frac{1}{2\left|\Omega\right|}.\label{eq: normal velocity}
\end{equation}
Later we will show how to use this velocity $V_{n}$ to find the optimal
domain which maximizes normalized $k-$th Steklov eigenvalue with
respect to the area for a given $k.$

\section{Steklov Eigenvalue Problems on the Complex Plane}

\subsection{On a Simply-Connected Domain}

In this section, we formulate the Steklov eigenvalue problem on the
complex plane $\mathbb{C}$ instead of $\mathbb{\mathbb{R}}^{2}.$
Consider the Steklov eigenvalue problem (\ref{eq:Steklov eigenvalue})
on a simply-connected domain $\Omega\subset\mathbb{C}.$ Due to the
Riemann Mapping Theorem that guarantees the existence of a unique
conformal mapping between any two simply-connected domains, we denote
$f=f(\omega)$ as the mapping function that maps the interior of a
unit circle $|\omega|=1$ where $w=re^{i\theta}=\xi+i\eta$ to the
interior of $\Omega$. Furthermore, every harmonic function is the
real part of an analytic function, $u=\Re\{\Psi\}$ where $\Psi$
is the complex potential and $\Re\{\Psi\}$ denotes the real part
of the argument $\Psi$. The advantage of this formulation is that
we no longer need to solve the equation on $\Omega$ as $u$ satisfies
the Laplace's equation automatically. We only need to find the solution
satisfies the boundary condition. 
\begin{figure}
\begin{centering}
\includegraphics[scale=0.8]{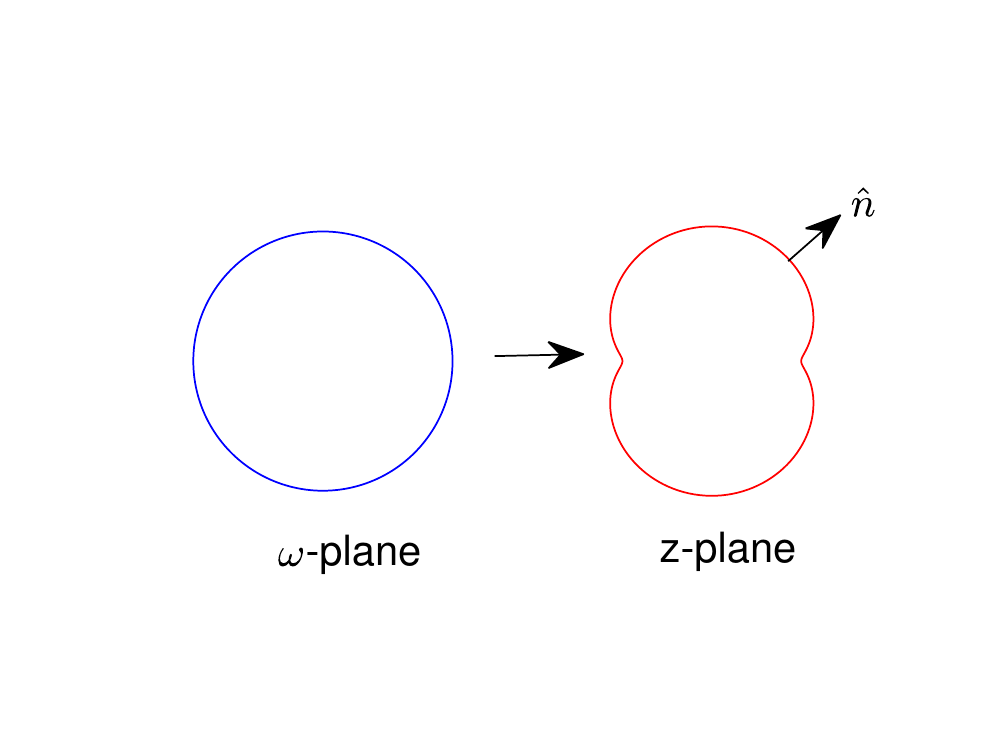}
\par\end{centering}
\caption{The mapping from a unit circle on $\omega-$plane to a simply-connected
domain on $z-$plane.\label{fig:The-mapping-from w to z}}
\end{figure}

Parametrizing the boundary of the original domain $\Omega$ with $z(\theta)=x(\theta)+iy(\theta)=f(\omega)$,
$\mid\omega\mid=1$ as shown in Figure \ref{fig:The-mapping-from w to z}.
The outward unit normal is 
\[
\hat{n}=(\frac{\dot{y}}{\sqrt{\dot{y}^{2}+\dot{x}^{2}}},\frac{-\dot{x}}{\sqrt{\dot{y}^{2}+\dot{x}^{2}}})
\]
where $\dot{x}=\frac{dx}{d\theta}$, $\dot{y}=\frac{dy}{d\theta}$
and the gradient of $u$ is 
\[
\nabla_{z}u=u_{x}+iu_{y}.
\]
Thus the derivative in the normal direction is given by 
\begin{equation}
\hat{n}\cdot\nabla_{z}u=\Im\left\{ (\frac{\dot{x}+i\dot{y}}{\sqrt{\dot{x}^{2}+\dot{y}^{2}}})(u_{x}-iu_{y})\right\} =\Im\left\{ \frac{\dot{z}}{\left|f_{\omega}\right|}\overline{\nabla_{z}u}\right\} \label{eq:normal derivatives}
\end{equation}
where $\Im(\cdot)$ denotes the imaginary part of the argument. Since,
$z=f(w)$, we have $\dot{z}=f_{\omega}\dot{\omega}=if_{\omega}\omega$
and $\overline{\nabla_{z}u}=\Psi_{z}=\Psi_{\omega}/f_{\omega}.$ Thus,
we get 
\[
\hat{n}\cdot\nabla_{z}u=\Im\left\{ \frac{if_{\omega}\omega}{\left|f_{\omega}\right|}\frac{\Psi_{\omega}}{f_{\omega}}\right\} =\Re\left\{ \frac{\omega}{\left|f_{\omega}\right|}\Psi_{\omega}\right\} \quad\text{on\ensuremath{\,\left|\omega\right|=1}}.
\]
The boundary condition $\frac{\partial u}{\partial n}=\lambda u$
in (\ref{eq:Steklov eigenvalue}) thus becomes
\begin{equation}
\Re\left\{ \omega\Psi_{\omega}\right\} =\lambda\left|f_{\omega}\right|\Re\left\{ \Psi\right\} \qquad\text{on}\left|\omega\right|=1.\label{eq: bc eigenvalue}
\end{equation}
Note that $\lambda=0$ is an eigenvalue and its corresponding eigenfunction
$u=\Re\{\Psi\}$ is a constant function. In this formulation, it is
not necessary to solve the harmonic equation as the real part of an
analytic function is always harmonic. However, it is required to know
the mapping function $f(\omega)$ and solve the equation (\ref{eq: bc eigenvalue})
on the unit circle. \textcolor{black}{In some cases, it is not easy
to find a conformal mapping between an arbitrary simply-connected
domain and the unit circle. When this happens, the Schwarz\textendash Christoffel
transformation \cite{driscoll1994schwarz} can be used to estimate
the mapping. }

\subsection{Steklov Eigenvalues of an Annulus}

In Section \ref{subsec:On-an-annulus}, we find Steklov eigenvalues
on an annulus $\Omega=B(0,1)\setminus B(0,\epsilon)$ in $\mathbb{R}^{2}$.
Here we reformulate the same problem in $\mathbb{C}$ and show that
the same equation is obtained for determining the eigenvalues. The
boundary conditions (\ref{eq:ellipse eigenvalue equation}) in the
complex formula are
\begin{eqnarray}
\Re(\omega\Psi_{\omega}) & = & \,\,\,\lambda\Re(\Psi),\;\text{on \ensuremath{\left|\omega\right|=1},}\nonumber \\
\Re(\omega\Psi_{\omega}) & = & -\lambda\Re(\Psi),\;\text{on}\ensuremath{\left|\omega\right|=\epsilon.}\label{eq: ellipse equation on complex plane}
\end{eqnarray}
where $\omega=re^{i\theta}$. Plugging $\Psi=\underset{k}{\sum}a_{k}\omega^{k}$
into (\ref{eq: ellipse equation on complex plane}) leads to
\[
\sum_{k}k(a_{k}-\overline{a_{-k}})e^{ik\theta}=\lambda\sum_{k}(a_{k}+\overline{a_{-k}})e^{ik\theta},
\]

\[
\sum_{k}k(a_{k}\epsilon^{k-1}-\overline{a_{-k}}\epsilon^{-k-1})e^{ik\theta}=-\lambda\sum_{k}(a_{k}\epsilon^{k}+\overline{a_{-k}}\epsilon^{-k})e^{ik\theta},
\]
which implies that

\[
\left[\begin{array}{cc}
\lambda\epsilon^{k}+k\epsilon^{k-1} & \lambda\epsilon^{-k}-k\epsilon^{-k-1}\\
\lambda-k & \lambda+k
\end{array}\right]\left[\begin{array}{c}
a_{k}\\
\overline{a_{-k}}
\end{array}\right]=\left[\begin{array}{c}
0\\
0
\end{array}\right].
\]
As shown in Section \ref{subsec:On-an-annulus}, the Steklov eigenvalues
are determined by the roots of the polynomial (\ref{eq:poly eqn}).

\subsection{Shape Optimization Problem}

Here we formulate the velocity (\ref{eq: normal velocity}) in the
complex plane $\mathbb{C}.$ By using the fact that $\left|\omega\right|^{2}=\bar{\omega}\omega=1$,
we obtain
\[
\frac{d}{dt}\left(\bar{\omega}\omega\right)=\bar{\omega}_{t}\omega+\bar{\omega}\omega_{t}=2\Re\{\bar{\omega}_{t}\omega\}=0.
\]
Since the mapping function is $z=f(\omega,t)$, we have $\frac{dz}{dt}=f_{t}+f_{\omega}\omega_{t}.$
The normal component of the velocity is given by
\begin{align*}
V_{n} & =\hat{n}\cdot V=\Im\left\{ \frac{if_{\omega}\omega}{|f_{\omega}|}\overline{\left(f_{\omega}\omega_{t}+f_{t}\right)}\right\} =\Re\left\{ \frac{f_{\omega}\omega}{|f_{\omega}|}\overline{\left(f_{\omega}\omega_{t}+f_{t}\right)}\right\} \\
 & =\Re\left\{ \frac{\overline{f_{\omega}\omega}}{|f_{\omega}|}\left(f_{\omega}\omega_{t}+f_{t}\right)\right\} =\Re\left\{ |f_{\omega}|\overline{\omega}\omega_{t}+\frac{f_{t}\left|f_{\omega}\right|}{\omega f_{\omega}}\right\} \\
 & =\Re\left\{ \frac{f_{t}\left|f_{\omega}\right|}{\omega f_{\omega}}\right\} .
\end{align*}
Therefore, the velocity (\ref{eq: normal velocity}) becomes

\[
\Re\left\{ \frac{f_{t}\mid f_{w}\mid}{wf_{w}}\right\} =\mid\nabla u\mid^{2}-2\lambda^{2}u^{2}-\lambda\kappa u^{2}+\frac{\lambda}{2|\Omega|}
\]
where $u$ is the normalized eigenfunction satisfying 
\begin{equation}
\int_{\partial\varOmega}u^{2}ds=\int_{|\omega|=1}\left(\Re\left\{ \Psi\right\} \right)^{2}|f_{\omega}|d\omega=1.\label{eq: normalized eigenfunction}
\end{equation}
Thus

\begin{equation}
\Re\left\{ \frac{f_{t}}{wf_{w}}\right\} =R(f,\Psi)\label{eq:moving bc}
\end{equation}
where the right hand side function $R(f,\Psi)$ is 
\[
R(f,\Psi)=\frac{1}{\left|f_{\omega}\right|}\left(\mid\Psi_{\omega}\mid^{2}\frac{1}{\left|f_{\omega}\right|^{2}}-2\lambda^{2}\Re(\Psi)^{2}-\lambda\kappa\Re(\Psi)^{2}+\frac{\lambda}{2|\Omega|}\right),
\]
$|\Omega|$ is the area of the given domain and the curvature is 
\begin{equation}
\kappa=\frac{\Re\left\{ \overline{\omega f_{\omega}}\left(\omega\left(\omega f_{\omega}\right)_{\omega}\right)\right\} }{\left|\omega f_{\omega}\right|^{3}}.\label{eq: curvature}
\end{equation}

Now, since $f$ is analytic in $\mid\omega\mid<1$ , 
\[
\frac{f_{t}}{\omega f_{\omega}}\text{ is analytic in}\left|\omega\right|<1.
\]
By using the Poisson integral formula, the value of an analytic function
in the domain $\mid\omega\mid<1$ can be obtained in term of its real
part evaluated on the unit circle. The equation (\ref{eq:moving bc})
implies that 
\begin{align*}
\frac{f_{t}}{\omega f_{\omega}} & =\frac{1}{2\pi i}\oint_{\mid\omega^{\prime}\mid=1}\frac{1}{\omega^{\prime}}\frac{\omega^{\prime}+\omega}{\omega\prime-\omega}\Re\left\{ R\left(f\left(\omega'\right),\Psi\left(\omega'\right)\right)\right\} d\omega^{\prime}\\
 & =\Re\left\{ R\left(f\left(\omega\right),\Psi\left(\omega\right)\right)\right\} +i\text{\ensuremath{\mathcal{H}}}\left\{ R\left(f\left(\omega\right),\Psi\left(\omega\right)\right)\right\} 
\end{align*}
where 
\[
\text{\ensuremath{\mathcal{H}}}[R\left(f\left(e^{i\theta}\right),\Psi\left(e^{i\theta}\right)\right)]=\frac{1}{2\pi}\int_{-\pi}^{\pi}\cot(\frac{\theta^{\prime}-\theta}{2})\Re\left\{ R\left(f\left(e^{i\theta'}\right),\Psi\left(e^{i\theta'}\right)\right)\right\} d\theta^{\prime}
\]
is the Hilbert transform. Thus we have 
\begin{equation}
f_{t}=\omega f_{\omega}(\Re(R\left(f\left(\omega\right),\Psi\left(\omega\right)\right))+i\text{\ensuremath{\mathcal{H}}}[R\left(f\left(\omega\right),\Psi\left(\omega\right)\right)])\label{eq: dynamic equation}
\end{equation}
which provides the deformation of the domain via the changes of the
conformal mapping. 

\section{Numerical Approaches for Solving Steklov Eigenvalue Problems}

In this section, we discuss the details of numerical discretization.
Assume $f$ and $\psi$ are represented as series expansions, i.e.
\[
f(w)=\sum_{-\infty}^{\infty}a_{k}\omega^{k}\text{ and }\Psi=\sum_{-\infty}^{\infty}c_{k}\omega^{k},
\]
respectively. In Section \ref{subsec:Forward-Solvers}, we discuss
how to find Steklov eigenvalues and eigenfunctions on a given domain
which is represented by $z=f(\omega)$, $|\omega|\le1$. This requires
to find eigenvalues $\lambda$ and analytic functions $\Psi$ whose
real part are eigenfunctions in Equation (\ref{eq: bc eigenvalue})
for a given $f.$ In Section \ref{subsec:Optimization-Solvers}, we
discuss how to discretize Equation (\ref{eq:moving bc}) on a unit
circle to obtain a system of ordinary differential equations (ODEs)
of the coefficients $a_{k}(t)$ of $f(\omega,t)$ with a given initial
guess of $a_{k}(0)$ of $f(\omega,0)$. The stationary state of this
system of ODEs gives the optimal area-normalized Steklov eigenvalue. 

\subsection{Forward Solvers\label{subsec:Forward-Solvers}}

\textcolor{black}{Given }$f(w)=\sum_{-\infty}^{\infty}a_{k}w^{k}$,
we solve (\ref{eq: bc eigenvalue}) numerically on $|\omega|=1$ by
parametrizing the unit circle by using the angle $\theta$
\[
\omega=e^{i\theta},\:\theta=[0,2\pi).
\]
Note that $a_{k}=0$ for $k<0$ as the domain is mapping to the interior
of the unit circle, i.e. $|\omega|\le1.$ The derivative of $f$ can
be obtained as 
\[
f_{\omega}=\sum_{-\infty}^{\infty}ka_{k}\omega^{k-1}
\]
and the magnitude of $|f_{\omega}|=\left(f_{\omega}\overline{f_{\omega}}\right)^{\frac{1}{2}}$
can be obtained in a series expansion again. Assume that the series
expansion of $|f_{\omega}|$ is 
\[
|f_{\omega}|=\sum_{-\infty}^{\infty}d_{l}\omega^{l}.
\]
Since $|f_{\omega}|$ is real, we must have $d_{l}=\overline{d_{-l}}$.
Denote the expansion of $\Psi$ as 
\[
\Psi=\sum_{-\infty}^{\infty}c_{k}\omega^{k}
\]
where $c_{k}=0$ for $k<0$ too. Plugging these series expansions
into (\ref{eq: bc eigenvalue}), we have 
\[
\sum_{-\infty}^{\infty}k(c_{k}-\overline{c_{-k}})\omega^{k}=\lambda\sum_{m=-\infty}^{\infty}\sum_{l=-\infty}^{\infty}(c_{m}+\overline{c_{-m}})d_{l}\omega^{m+l}.
\]
We can then use the identity 
\[
\sum_{m=-\infty}^{\infty}\sum_{l=-\infty}^{\infty}(c_{m}+\overline{c_{-m}})d_{l}\omega^{m+l}=\sum_{k=-\infty}^{\infty}\sum_{m=-\infty}^{\infty}(c_{m}+\overline{c_{-m}})d_{k-m}\omega^{k}.
\]
By matching the coefficients of $\omega^{k}$, we have 
\begin{equation}
\lambda\sum_{m=-\infty}^{\infty}(c_{m}+\overline{c_{-m}})d_{k-m}=k(c_{k}-\overline{c_{-k}}).\label{eq: continuous conformall FF}
\end{equation}
Denote the real and complex part of $c_{n}$, $d_{n}$ by $c_{n}^{r}$,
$d_{n}^{r}$ and $c_{n}^{i}$, $d_{n}^{i}$, respectively, we then
have 
\[
\lambda\sum_{m=-\infty}^{\infty}(c_{m}^{r}+ic_{m}^{i}+c_{-m}^{r}-ic_{-m}^{i})(d_{k-m}^{r}+id_{k-m}^{i})=k(c_{k}^{r}+ic_{k}^{i}-c_{-k}^{r}+ic_{-k}^{i}).
\]
By comparing real and imaginary parts, we have 
\begin{equation}
\left\{ \begin{array}{ccc}
\lambda\sum_{m=-\infty}^{\infty}c_{m}^{r}(d_{k-m}^{r}+d_{k+m}^{r})+c_{m}^{i}(d_{k+m}^{i}-d_{k-m}^{i}) & = & kc_{k}^{r},\\
\lambda\sum_{m=-\infty}^{\infty}c_{m}^{r}(d_{k-m}^{i}+d_{k+m}^{i})+c_{m}^{i}(d_{k-m}^{r}-d_{k+m}^{r}) & = & kc_{k}^{i}.
\end{array}\right.\label{eq: system of c}
\end{equation}
In numerical computation, the series expansion is carried out numerically
by truncating the series expansion at $k=\frac{N}{2}$ and Fast Fourier
Transform (FFT) is used to efficiently compute quantities in $\omega-$plane
and $z-$plane. Denote $\frac{N}{2}$ as $N_{2}$. Thus
\[
f(\omega)\approx\sum_{k=-N_{2}}^{N_{2}}a_{k}\omega^{k}=\sum_{k=-N_{2}}^{N_{2}}a_{k}e^{ik\theta}
\]
and
\[
f_{\omega}=\sum_{k=-N_{2}}^{N_{2}}ka_{k}\omega^{k-1}.
\]
Denote 
\[
|f_{\omega}|=\sum_{l=-N}^{N}d_{l}\omega^{l}
\]
where $d_{l}$, $-N\le l\le N$, are obtained by using the pseudo-spectral
method. We use inverse Fourier transform (IFFT) to obtain $f_{\omega}$
in physical space and compute $|f_{\omega}|$ in physical space, then
use FFT to get $d_{l}$ in Fourier space. The aliasing of a nonlinear
produ\textcolor{black}{ct is avoided by adopting the zero-padding. }

The system of infinite equations (\ref{eq: system of c}) is approximated
by the system of finite equations for\textcolor{red}{{} }\textcolor{black}{0:$N_{2}$-modes
w}hich gives
\begin{equation}
\lambda AC=BC\label{eq:system}
\end{equation}
where
\begin{align*}
A_{k+1,m+1} & =d_{k-m}^{r}+d_{k+m}^{r},\quad\text{for}\quad0\le k\le N_{2},0\le m\le N_{2},\\
A_{k+1,m+N_{2}+1} & =-d_{k-m}^{i}+d_{k+m}^{i},\quad\text{for}\quad0\le k\le N_{2},1\le m\le N_{2},\\
A_{k+N_{2}+1,m+1} & =d_{k-m}^{i}+d_{k+m}^{i},\quad\text{for}\quad1\le k\le N_{2},0\le m\le N_{2},\\
A_{k+N_{2}+1,m+N_{2}+1} & =d_{k-m}^{r}-d_{k+m}^{r},\quad\text{for}\quad1\le k\le N_{2},1\le m\le N_{2},
\end{align*}
and
\begin{align*}
B_{k+1,m+1} & =k\delta_{k,m},\quad\text{for}\quad0\le k\le N_{2},0\le m\le N_{2},\\
B_{k+1,m+N_{2}+1} & =0,\quad\text{for}\quad0\le k\le N_{2},1\le m\le N_{2},\\
B_{k+N_{2}+1,m+1} & =0,\quad\text{for}\quad1\le k\le N_{2},0\le m\le N_{2},\\
B_{k+N_{2}+1,m+N_{2}+1} & =k\delta_{k,m},\quad\text{for}\quad1\le k\le N_{2},1\le m\le N_{2},
\end{align*}
and 
\[
C=\left[\begin{array}{c}
C^{r}\\
C^{i}
\end{array}\right]
\]
where 
\[
C^{r}=\left[\begin{array}{c}
c_{0}^{r}\\
c_{1}^{r}\\
c_{2}^{r}\\
\vdots\\
c_{N_{2}}^{r}
\end{array}\right],\quad C^{i}=\left[\begin{array}{c}
c_{0}^{i}\\
c_{1}^{i}\\
c_{2}^{i}\\
\vdots\\
c_{N_{2}}^{i}
\end{array}\right].
\]
By solving the linear system (\ref{eq:system}) we could find the
coefficient vector $C$ and its corresponding eigenvalue $\lambda$.
We assign zero values for $c_{k}^{r}$ and $c_{k}^{i}$ for $k>N_{2}.$
Thus the corresponding eigenfunction will be given by $u=\Re\{\Psi\}=\Re\{\sum_{-N_{2}}^{N_{2}}c_{k}\omega^{k}\}.$

Now, if we assume that the coefficients $d_{n}$ are real we will
be able to reduce the matrix size and solve the problem even more
efficiently. In this case, we have 

\begin{equation}
\left\{ \begin{array}{ccc}
\lambda\sum_{m=0}^{N_{2}}c_{m}^{r}(d_{k-m}+d_{k+m}) & = & kc_{k}^{r},\\
\lambda\sum_{m=1}^{N_{2}}c_{m}^{i}(d_{k-m}-d_{k+m}) & = & kc_{k}^{i}.
\end{array}\right.\label{eq: Discrete system 2}
\end{equation}
The 0:$N_{2}$-modes approximation gives
\[
\lambda A^{r}C^{r}=B^{r}C^{r},\lambda A^{i}C^{i}=B^{i}C^{i},
\]
where 
\[
A_{k+1,m+1}^{r}=d_{k-m}+d_{k+m},\quad B_{k+1,m+1}^{r}=k\delta_{k,m}\quad\text{for}\quad0\le k\le N_{2},0\le m\le N_{2},
\]

\[
A_{k,m}^{i}=d_{k-m}-d_{k+m},\quad B_{k,m}^{i}=k\delta_{k,m},\quad\text{for}\quad0\le k\le N_{2},0\le m\le N_{2},
\]
and
\[
C^{r}=\left[\begin{array}{c}
c_{0}^{r}\\
c_{1}^{r}\\
c_{2}^{r}\\
\vdots\\
c_{N_{2}}^{r}
\end{array}\right],\quad C^{i}=\left[\begin{array}{c}
c_{0}^{i}\\
c_{1}^{i}\\
c_{2}^{i}\\
\vdots\\
c_{N_{2}}^{i}
\end{array}\right].
\]

\subsection{Optimization Solvers\label{subsec:Optimization-Solvers}}

In this section, we discuss how to solve the dynamic equation (\ref{eq: dynamic equation})
by method of lines and spectral method in the variable $\omega$.
Given a conformal mapping $f(\omega,t)=\sum_{-N_{2}}^{N_{2}}a_{k}(t)\omega^{k}$,
we use the method discussed in \ref{subsec:Forward-Solvers} to obtain
$k$th eigenvalue $\lambda_{k}$, its corresponding eigenfunction
$u_{k}=\Re\{\Psi\}$ where $\Psi(w,t)=\sum_{-N_{2}}^{N_{2}}c_{k}(t)\omega^{k}$.
Notice that this eigenfunction is not normalized. To find the normalization
constant, we compute the Fourier coefficient representation of 
\[
\left(\Re\left\{ \Psi\right\} \right)^{2}|f_{\omega}|=\sum_{-N_{2}}^{N_{2}}b_{k}(t)\omega^{k}
\]
via a pseudo-spectral method and then the normalization condition
(\ref{eq: normalized eigenfunction}) is approximated by 
\[
\int_{|\omega|=1}\left(\Re\left\{ \Psi\right\} \right)^{2}|f_{\omega}|d\omega\approx2\pi b_{0}(t).
\]
The normalized eigenfunction $u=\Re\left\{ \tilde{\Psi}\right\} $
where 
\[
\tilde{\Psi}(w,t)=\sum_{-N_{2}}^{N_{2}}\frac{1}{\sqrt{2\pi b_{0}}}c_{k}(t)\omega^{k}=\sum_{-N_{2}}^{N_{2}}\tilde{c}_{k}(t)\omega^{k}.
\]
The curvature term can be computed via the formula (\ref{eq: curvature})
by using the following expansions 
\[
\omega f_{\omega}=\sum_{-N_{2}}^{N_{2}}ka_{k}(t)\omega^{k}
\]
\[
\omega\left(\omega f_{\omega}\right)_{\omega}=\sum_{-N_{2}}^{N_{2}}k^{2}a_{k}(t)\omega^{k}.
\]
The area term is obtained by 
\[
|\Omega|=\sum_{-N_{2}}^{N_{2}}\pi k|a_{k}|^{2}.
\]
Plugging 
\[
|f_{\omega}|=\sum_{-N}^{N}d_{l}\omega^{l},\quad\tilde{\Psi}_{\omega}=\sum_{-N_{2}}^{N_{2}}k\tilde{c}_{k}(t)\omega^{k},
\]
the eigenvalue, the curvature, and the area into the right hand side
of (\ref{eq:moving bc}), we obtain $R(f,\tilde{\Psi})$ in terms
of Fourier series. All the nonlinear term is obtained by using pseudo-spectral
method. We then use discrete Hilbert transform to find the complex
conjugate of $R(f,\tilde{\Psi})$ and then compute the right hand
side of (\ref{eq: dynamic equation}). Denote the series expansion
of the right hand side as 
\[
wf_{w}(\Re\left\{ R\left(f(\omega),\tilde{\Psi}(\omega\right)\right\} )+i\text{\ensuremath{\mathcal{H}}}[R\left(f(\omega),\tilde{\Psi}(\omega\right)])=\sum_{-N_{2}}^{N_{2}}r_{k}(t)\omega^{k}.
\]
Note that $r_{k}$ depends on time and $a_{k}$, $-N_{2}\le k\le N_{2}$.
Since $f_{t}(\omega,t)=\sum_{-\infty}^{\infty}a_{k}^{\prime}(t)\omega^{k}$,
the dynamic equation (\ref{eq: dynamic equation}) becomes a system
of $N+1$ nonlinear ODEs in Fourier Coefficients

\begin{equation}
a_{k}^{\prime}(t)=r_{k}(t),-N_{2}\le k\le N_{2}.\label{eq: a_k ODE}
\end{equation}

\section{Numerical Results}

\subsection{Forward Solvers}

Here we first test our forward solvers on various domains to demonstrate
the spectral convergence of the numerical approaches described in
Section \ref{subsec:Forward-Solvers}. We verify the accuracy of the
code by testing the first $12$ eigenvalues on smooth shapes. 

\subsubsection{Steklov Eigenvalues on a Unit Disk }

When we consider the unit circle, the mapping function is $f(\omega)=\omega$
which gives $|f_{\omega}|=1.$ Thus $d_{0}=1$ and $d_{l}=0$ for
all $l\neq0.$ The system of equations (\ref{eq: Discrete system 2})
becomes 

\begin{equation}
\left\{ \begin{array}{ccc}
\lambda c_{k}^{r}=kc_{k}^{r}, &  & k=0,1,2,3,...\\
\lambda c_{k}^{i}=kc_{k}^{i}, &  & k=1,2,3,....
\end{array}\right.\label{eq: Discrete system 2-1}
\end{equation}
If $\lambda=0$, $c_{k}^{r}=c_{k}^{i}=0$ for all positive integer
and $c_{0}^{r}$ is an arbitrary constant. If $\lambda$ is a particular
integer $k_{1}$, i.e., $\lambda=k_{1},$ we must have 
\[
c_{k}^{r}=c_{k}^{i}=0\quad\text{for}\quad k\neq k_{1},
\]
and $c_{k_{1}}^{r}$ and $c_{k_{1}}^{i}$ are arbitrary constants.
Thus, Steklov eigenvalue for the unit circle are 
\[
0,1,1,2,2,3,3,\cdots,k_{1},k_{1},\cdots.
\]
The normalized eigenvalues $\lambda_{k}^{A}(\Omega)=\lambda\sqrt{\mid\Omega\mid}$
are listed in Table \ref{tab:unit circle}. It is clear that spectral
accuracy is observed from the numerical results and the errors only
contain round off errors $O(10^{-16})$ on double-precision arithmetic. 

\begin{table}[H]
\begin{centering}
\begin{tabular}{|c|c|c|c|c|}
\hline 
$N$ & 2$^{4}$ & 2$^{5}$ & 2$^{12}$ & Exact\tabularnewline
\hline 
\hline 
$\lambda_{0}$ & 0 & 0 & 0 & $0$\tabularnewline
\hline 
$\lambda_{1}$ & 1.772453850905515 & 1.772453850905515 & 1.772453850905515 & 1.772453850905516\tabularnewline
\hline 
$\lambda_{2}$ & 1.772453850905515 & 1.772453850905515 & 1.772453850905515 & 1.772453850905516\tabularnewline
\hline 
$\lambda_{3}$ & 3.544907701811031 & 3.544907701811031 & 3.544907701811031 & 3.544907701811032\tabularnewline
\hline 
$\lambda_{4}$ & 3.544907701811031 & 3.544907701811031 & 3.544907701811031 & 3.544907701811032\tabularnewline
\hline 
$\lambda_{5}$ & 5.317361552716547  & 5.317361552716547  & 5.317361552716547  & 5.317361552716548\tabularnewline
\hline 
$\lambda_{6}$ & 5.317361552716547  & 5.317361552716547  & 5.317361552716547  & 5.317361552716548\tabularnewline
\hline 
$\lambda_{7}$ & 7.089815403622062 & 7.089815403622062 & 7.089815403622062 & 7.089815403622064\tabularnewline
\hline 
$\lambda_{8}$ & 7.089815403622062 & 7.089815403622062 & 7.089815403622062 & 7.089815403622064\tabularnewline
\hline 
$\lambda_{\text{9}}$ & 8.862269254527577 & 8.862269254527577 & 8.862269254527577 & 8.862269254527579\tabularnewline
\hline 
$\lambda_{10}$ & 8.862269254527577 & 8.862269254527577 & 8.862269254527577 & 8.862269254527579\tabularnewline
\hline 
$\lambda_{11}$ & 10.634723105433094 & 10.634723105433094 & 10.634723105433094 & 10.634723105433096\tabularnewline
\hline 
\end{tabular}
\par\end{centering}
\caption{The first 12 eigenvalues $\lambda_{k},k=0,\ldots,11$  for different
numbers of grid points $N=2^{n},n=4,5,12$ on a unit circle.\label{tab:unit circle}}
\end{table}

\subsubsection{Steklov Eigenvalues on a Shape with 2-Fold Rotational Symmetry}

\begin{figure}[H]
\begin{centering}
\subfloat[]{\begin{centering}
\includegraphics[scale=0.6]{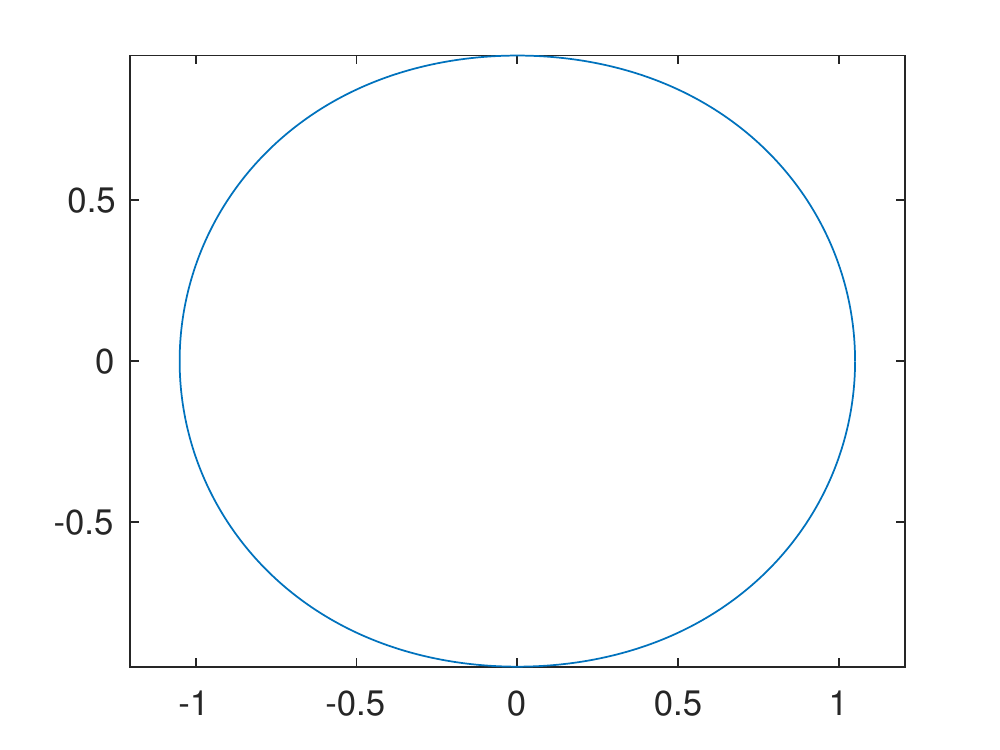}
\par\end{centering}
}\subfloat[]{\begin{centering}
\includegraphics[scale=0.6]{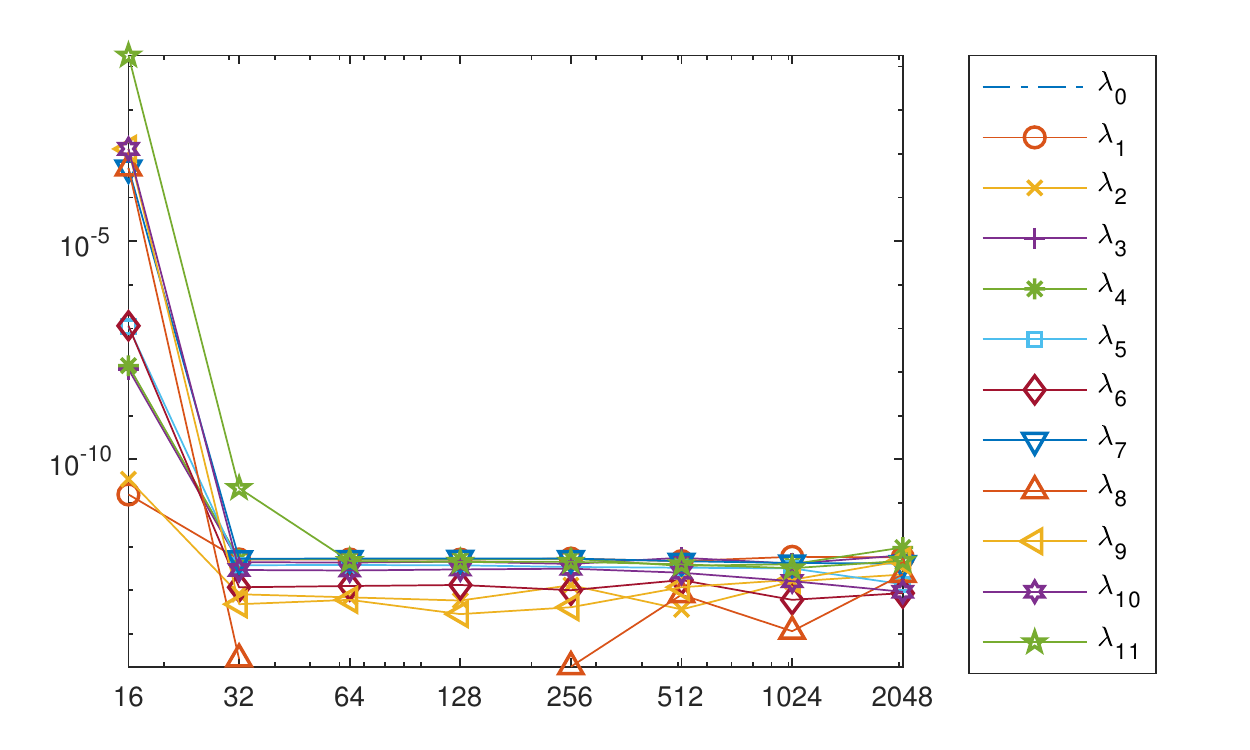}
\par\end{centering}
}
\par\end{centering}
\caption{(a) The 2-fold rotational symmetry shape with\textcolor{magenta}{{}
}$f(w)=w+0.05w^{3}$,$|\omega|\le1.$ (b) The log-log plot of errors
for the first 11 non-zero eigenvalues versus number of grid points
$N=2^{n},n=4,\ldots,11.$ \label{fig:elipse}}
\end{figure}

We use the mapping $f(w)=w+0.05w^{3}$ to generate a shape with 2-fold
rotational symmetry as shown in Figure \ref{fig:elipse}(a). In Table
\ref{tab:ellipse} we summarize the numerical results of Steklov eigenvalues.
We use the eigenvalues computed by using $2^{12}$ grids as true eigenvalues
and show the log-log plot of errors of the first 12 eigenvalues, i.e.
\[
\text{error=}|\lambda_{k}^{N}-\lambda_{k}^{2^{12}}|,k=0,\ldots,11,
\]
versus number of grid points $N=2^{4},2^{5},...2^{11}$ in Figure
\ref{fig:elipse}(b). It is clear that the spectral accuracy is achieved. 

\begin{table}[H]
\begin{centering}
\begin{tabular}{|c|c|c|c|c|}
\hline 
$N$ & 2$^{4}$ & 2$^{5}$ & 2$^{6}$ & 2$^{7}$\tabularnewline
\hline 
\hline 
$\lambda_{0}$ & 0 & 0 & 0 & $0$\tabularnewline
\hline 
$\lambda_{1}$ & 1.643146123296456  & 1.643146123280263 & 1.643146123280263 & 1.643146123280268 \tabularnewline
\hline 
$\lambda_{2}$ & 1.904409864808107 & 1.904409864772927 & 1.904409864772939 & 1.904409864772950\tabularnewline
\hline 
$\lambda_{3}$ & 3.509482564053473 & 3.509482552385534 & 3.509482552385528 & 3.509482552385548\tabularnewline
\hline 
$\lambda_{4}$ & 3.567218990382545 & 3.567218976359059 & 3.567218976359065 & 3.567218976359050\tabularnewline
\hline 
$\lambda_{5}$ & 5.298764914769874 & 5.298764805372437 & 5.298764805372433 & 5.298764805372439\tabularnewline
\hline 
$\lambda_{6}$ & 5.316931803045312 & 5.316931688027542 & 5.316931688027550 & 5.316931688027557\tabularnewline
\hline 
$\lambda_{7}$ & 7.074710761837761 & 7.074238491011210 & 7.074238491011200 & 7.074238491011197\tabularnewline
\hline 
$\lambda_{8}$ & 7.079268312074488  & 7.078792636301956 & 7.078792636301953 & 7.078792636301953\tabularnewline
\hline 
$\lambda_{\text{9}}$ & 8.846269410836159 & 8.844970458352126 & 8.844970458352138 & 8.844970458352106 \tabularnewline
\hline 
$\lambda_{10}$ & 8.847598359495487 & 8.846297249970153 & 8.846297249970146 & 8.846297249970162\tabularnewline
\hline 
$\lambda_{11}$ & 10.793832137331764 & 10.614565359904542 & 10.614565359883139 & 10.614565359883118\tabularnewline
\hline 
\end{tabular}
\par\end{centering}
\begin{centering}
\begin{tabular}{|c|c|c|c|c|}
\hline 
$N$ & 2$^{8}$ & 2$^{9}$ & 2$^{10}$ & 2$^{12}$\tabularnewline
\hline 
\hline 
$\lambda_{0}$ & 0 & 0 & 0 & 0\tabularnewline
\hline 
$\lambda_{1}$ & 1.643146123296456 & 1.643146123280306 & 1.643146123280187 & 1.643146123280772\tabularnewline
\hline 
$\lambda_{2}$ & 1.904409864772878 & 1.904409864772972 & 1.904409864773167 & 1.904409864773008\tabularnewline
\hline 
$\lambda_{3}$ & 3.509482552385497 & 3.509482552385653 & 3.509482552385503 & 3.509482552385095\tabularnewline
\hline 
$\lambda_{4}$ & 3.567218976359074  & 3.567218976358907 & 3.567218976358941 & 3.567218976358544\tabularnewline
\hline 
$\lambda_{5}$ & 5.298764805372470  & 5.298764805372484 & 5.298764805372494 & 5.298764805372812\tabularnewline
\hline 
$\lambda_{6}$ & 5.316931688027525 & 5.316931688027596 & 5.316931688027485 & 5.316931688027425\tabularnewline
\hline 
$\lambda_{7}$ & 7.074238491011203 & 7.074238491011272 & 7.074238491011313 & 7.074238491011736\tabularnewline
\hline 
$\lambda_{8}$ & 7.078792636301955 & 7.078792636302032 & 7.078792636301965 & 7.078792636301953\tabularnewline
\hline 
$\lambda_{\text{9}}$ & 8.844970458352119 & 8.844970458352195 & 8.844970458352252 & 8.844970458352078\tabularnewline
\hline 
$\lambda_{10}$ & 8.846297249970174  & 8.846297249970114 & 8.846297249970023 & 8.846297249969862\tabularnewline
\hline 
$\lambda_{11}$ & 10.614565359883121 & 10.614565359883064  & 10.614565359882992 & 10.614565359882672\tabularnewline
\hline 
\end{tabular}
\par\end{centering}
\caption{The first 12 eigenvalues $\lambda_{k},k=0,\ldots,11$  for different
numbers of grid points $N=2^{n},n=4,\ldots,10,12$ on $f(w)=w+0.05w^{3}$,
$|\omega|\le1.$\label{tab:ellipse}}
\end{table}

\subsubsection{Steklov Eigenvalues on a Shape with 5-Fold Rotational Symmetry}

\begin{figure}[h]
\begin{centering}
\subfloat[]{\begin{centering}
\includegraphics[scale=0.6]{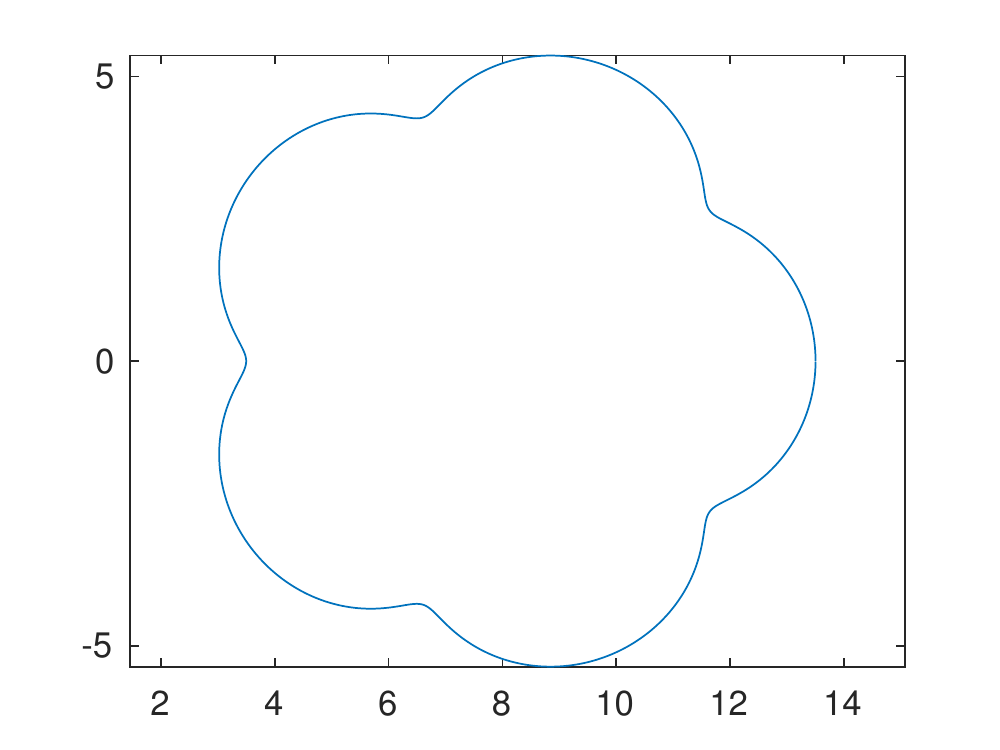}
\par\end{centering}

}\subfloat[]{\begin{centering}
\includegraphics[scale=0.6]{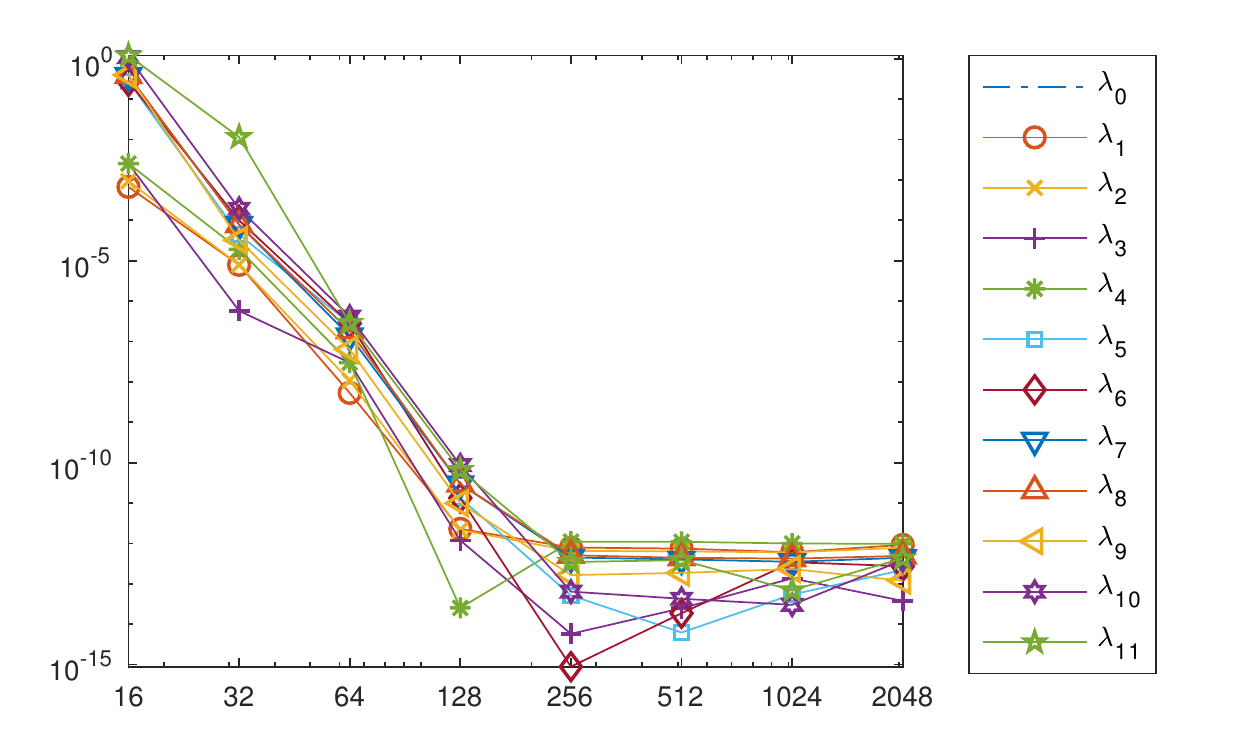}
\par\end{centering}

}
\par\end{centering}
\caption{(a) The 5-fold rotational symmetry shape with\textcolor{magenta}{{}
}$f(w)=8+5w+0.5w^{6}$,$|\omega|\le1.$ (b) The log-log plot of errors
for the first 11 non-zero eigenvalues versus number of grid points
$N=2^{n},n=4,\ldots,11.$ \label{fig:five-folds-}}
\end{figure}

We use the mapping $f(w)=8+5w+0.5w^{6}$ to generate a shape with
5-fold rotational symmetry as shown in Figure \ref{fig:five-folds-}(a).
In Table \ref{tab:5_folds}, we use the eigenvalues computed by using
$2^{12}$ grids as true eigenvalues and show the log-log plot of errors
of the first 12 eigenvalues, i.e.
\[
\text{error=}|\lambda_{k}^{N}-\lambda_{k}^{2^{12}}|,k=0,\ldots,11,
\]
versus number of grid points $N=2^{4},2^{5},...2^{11}$ in Figure
\ref{fig:five-folds-}(b). It is clear that the spectral accuracy
is achieved. 

\begin{table}[H]
\begin{centering}
\begin{tabular}{|c|c|c|c|c|}
\hline 
$N$ & 2$^{4}$ & 2$^{5}$ & 2$^{6}$ & 2$^{7}$\tabularnewline
\hline 
\hline 
$\lambda_{0}$ & 0 & 0 & 0 & $0$\tabularnewline
\hline 
$\lambda_{1}$ & 1.613981749710263 & 1.614659735134658 & 1.614651857980075 & 1.614651852652450\tabularnewline
\hline 
$\lambda_{2}$ & 1.615586942999712 & 1.614659740601958 & 1.614651863407194 & 1.614651852652469\tabularnewline
\hline 
$\lambda_{3}$ & 2.979901447266664 & 2.977376794062662 & 2.977377396867629 & 2.977377367030917\tabularnewline
\hline 
$\lambda_{4}$ & 2.979920850098075  & 2.977396410343160 & 2.977377396867634 & 2.977377367030926\tabularnewline
\hline 
$\lambda_{5}$ & 5.757902735512396 & 5.483423114699104  & 5.483379266795433 & 5.483378986137383\tabularnewline
\hline 
$\lambda_{6}$ & 5.757963817902539  & 5.483478476088597 & 5.483379266795448 & 5.483378986137454\tabularnewline
\hline 
$\lambda_{7}$ & 7.091240897150815 & 6.707817046860952 & 6.707738934321966 & 6.707738797445765\tabularnewline
\hline 
$\lambda_{8}$ & 7.092066936388594 & 6.707817092425547 & 6.707738981978962 & 6.707738797445767 \tabularnewline
\hline 
$\lambda_{\text{9}}$ & 8.053537400023426 & 7.657772022224528 & 7.657739872866688 & 7.657739809188358\tabularnewline
\hline 
$\lambda_{10}$ & 10.114031561463605 & 9.019776832943990 & 9.019583333936978 & 9.019582922824695\tabularnewline
\hline 
$\lambda_{11}$ & 11.339690354808871 & 10.150431507211664 & 10.138974110712084 & 10.138973824292398\tabularnewline
\hline 
\end{tabular}
\par\end{centering}
\begin{centering}
\begin{tabular}{|c|c|c|c|c|}
\hline 
$N$ & 2$^{8}$ & 2$^{9}$ & 2$^{10}$ & 2$^{12}$\tabularnewline
\hline 
\hline 
$\lambda_{0}$ & 0 & 0 & 0 & 0\tabularnewline
\hline 
$\lambda_{1}$ & 1.614651852650946 & 1.614651852650901 & 1.614651852650762 & 1.614651852650156\tabularnewline
\hline 
$\lambda_{2}$ & 1.614651852650962 & 1.614651852650941 & 1.614651852650909 & 1.614651852650308\tabularnewline
\hline 
$\lambda_{3}$ & 2.977377367029736 & 2.977377367029755 & 2.977377367029867 & 2.977377367029730\tabularnewline
\hline 
$\lambda_{4}$ & 2.977377367029792 & 2.977377367029804 & 2.977377367029905 & 2.977377367030901\tabularnewline
\hline 
$\lambda_{5}$ & 5.483378986124044 & 5.483378986123986  & 5.483378986124047 & 5.483378986123992\tabularnewline
\hline 
$\lambda_{6}$ & 5.483378986124095 & 5.483378986124115 & 5.483378986124439 & 5.483378986124096\tabularnewline
\hline 
$\lambda_{7}$ & 6.707738797416523 & 6.707738797416477 & 6.707738797416426 & 6.707738797416075\tabularnewline
\hline 
$\lambda_{8}$ & 6.707738797416656 & 6.707738797416588 & 6.707738797416567 & 6.707738797416147\tabularnewline
\hline 
$\lambda_{\text{9}}$ & 7.657739809178596 & 7.657739809178618 & 7.657739809178663 & 7.657739809178431\tabularnewline
\hline 
$\lambda_{10}$ & 9.019582922738280 & 9.019582922738174 & 9.019582922738246 & 9.019582922738216\tabularnewline
\hline 
$\lambda_{11}$ & 10.138973824227390 & 10.138973824227429 & 10.138973824227113 & 10.138973824227044\tabularnewline
\hline 
\end{tabular}
\par\end{centering}
\caption{The first 12 eigenvalues $\lambda_{k},k=0,\ldots,11$  for different
numbers of grid points $N=2^{n},n=4,\ldots,10,12.$\label{tab:5_folds}}
\end{table}

\subsubsection{Steklov Eigenvalues on a Cassini Oval.}

\begin{figure}[h]
\begin{centering}
\subfloat[]{\begin{centering}
\includegraphics[scale=0.6]{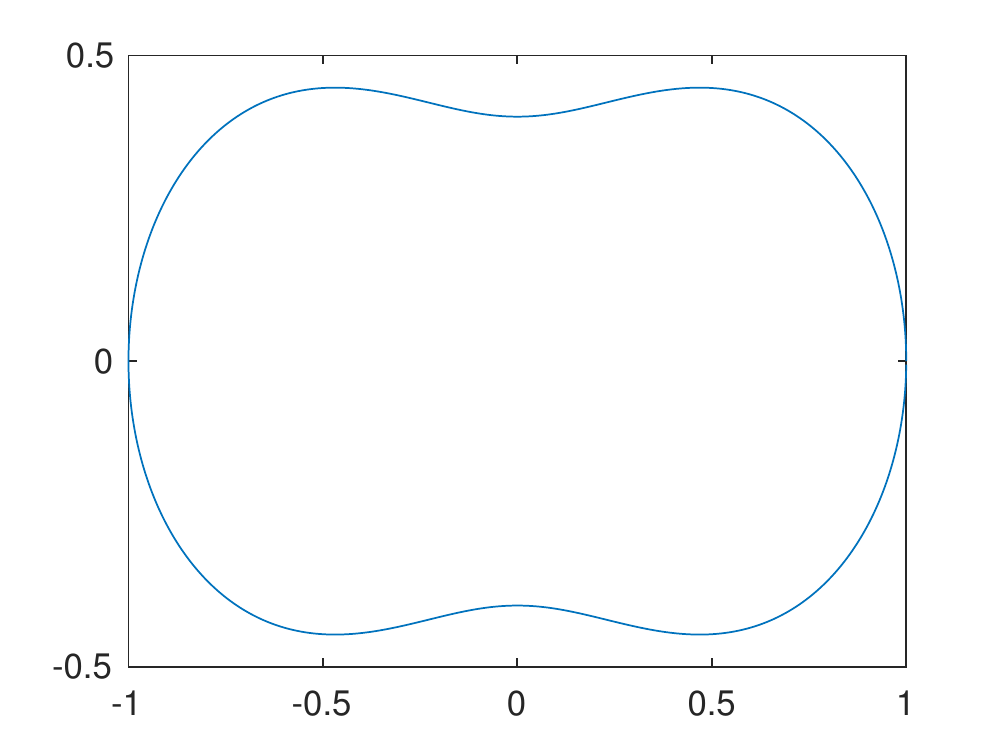}
\par\end{centering}
}\subfloat[]{\begin{centering}
\includegraphics[scale=0.6]{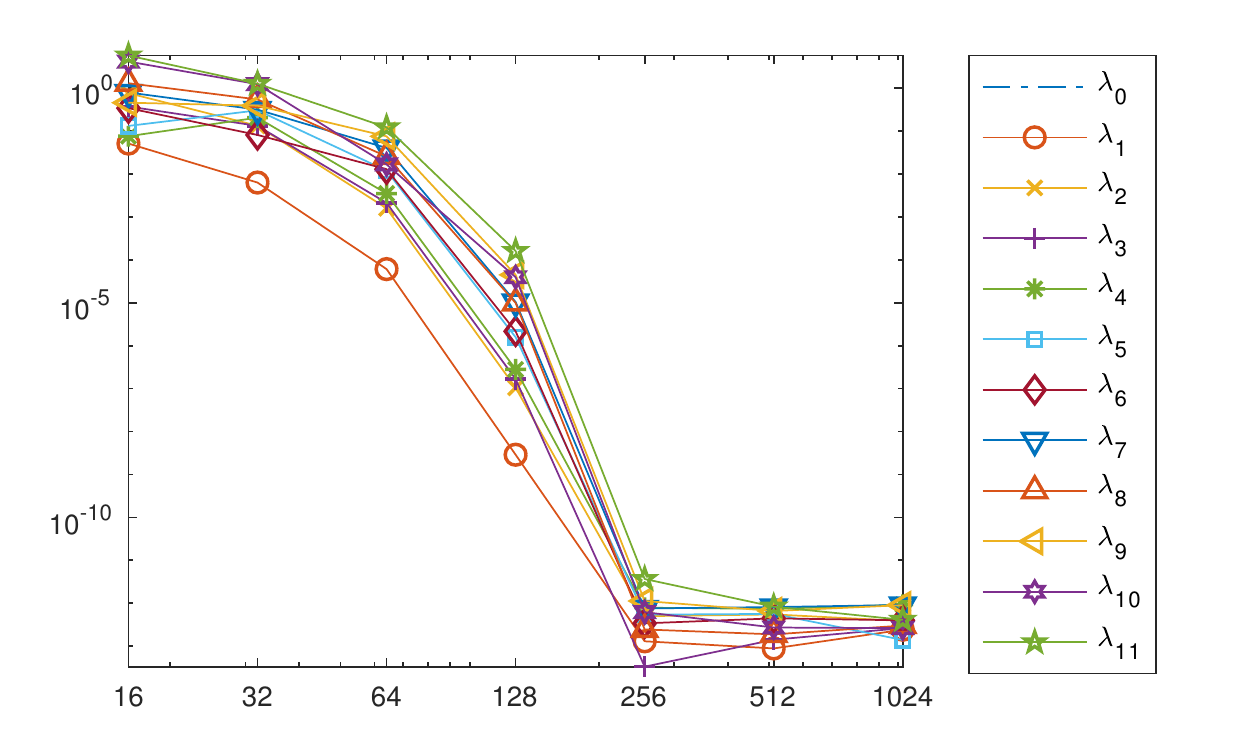}
\par\end{centering}
}
\par\end{centering}
\caption{(a) Cassini oval shape with $f(w)=\alpha w(\frac{2}{1+\alpha^{2}-(1-\alpha^{2})w^{2}})^{\frac{1}{2}}$
,$|\omega|\le1$, and $\alpha=0.4$. (b) The log-log plot of errors
for the first 11 non-zero eigenvalues versus number of grid points
$N=2^{n},n=4,\ldots,10.$ \label{fig:Cassini-Oval-shape}}
\end{figure}

All of aforementioned examples have finite terms expansion in $\omega$.
Here we show an example with infinite terms expansion in $\omega.$
The mapping $f(w)=\alpha w(\frac{2}{1+\alpha^{2}-(1-\alpha^{2})w^{2}})^{\frac{1}{2}}$
, where $\alpha=0.4$ is used to generate a Cassini Oval shape which
is shown in Figure \ref{fig:Cassini-Oval-shape}(a). In Table \ref{tab:cosaini oval}
we use the eigenvalues computed by using $2^{12}$ grids as true eigenvalues
and show the log-log plot of errors of the first 12 eigenvalues, i.e.
\[
\text{error=}|\lambda_{k}^{N}-\lambda_{k}^{2^{12}}|,k=0,\ldots,11,
\]
versus number of grid points $N=2^{4},2^{5},...2^{10}$ in Figure
\ref{fig:Cassini-Oval-shape}(b). It is also clear that the spectral
accuracy is achieved.

\begin{table}[H]
\begin{centering}
\begin{tabular}{|c|c|c|c|c|}
\hline 
$N$ & 2$^{4}$ & 2$^{5}$ & 2$^{6}$ & 2$^{7}$\tabularnewline
\hline 
\hline 
$\lambda_{0}$ & 0 & 0 & 0 & $0$\tabularnewline
\hline 
$\lambda_{1}$ & 0.872759997500228 & 0.827902995301854 & 0.821644770560566 & 0.821583902061334\tabularnewline
\hline 
$\lambda_{2}$ & 2.124401456784662 & 2.756054635303737 & 2.886951802792420 & 2.888537681079042\tabularnewline
\hline 
$\lambda_{3}$ & 2.571449696110635 & 3.077030643209814 & 2.946970106404462 & 2.944846781799040\tabularnewline
\hline 
$\lambda_{4}$ & 3.265821414464841 & 3.136596486946471 & 3.338218243465505 & 3.341726009279193\tabularnewline
\hline 
$\lambda_{5}$ & 4.418763601488365 & 4.854880686612021 & 4.562691155962767 & 4.550749526079698\tabularnewline
\hline 
$\lambda_{6}$ & 5.378320687602239 & 4.955570564610835 & 5.023787664833372 & 5.036737477441735\tabularnewline
\hline 
$\lambda_{7}$ & 7.025574559110670 & 6.548839953593787 & 6.273463980936640 & 6.233063933209499\tabularnewline
\hline 
$\lambda_{8}$ & 7.626882523375537 & 6.870625350063478 & 6.299773422917886 & 6.325481073833819\tabularnewline
\hline 
$\lambda_{\text{9}}$ & 8.263730860061949 & 8.196071839036918 & 7.881197306312033 & 7.805852388999546\tabularnewline
\hline 
$\lambda_{10}$ & 12.070450297339713 & 9.144520167848396  & 7.891680915128140 & 7.908376668589532\tabularnewline
\hline 
$\lambda_{11}$ & 15.149068049247919 & 10.668830636509751 & 9.526157620343742 & 9.404387869498899\tabularnewline
\hline 
\end{tabular}
\par\end{centering}
\begin{centering}
\begin{tabular}{|c|c|c|c|c|}
\hline 
$N$ & 2$^{8}$ & 2$^{9}$ & 2$^{10}$ & 2$^{12}$\tabularnewline
\hline 
\hline 
$\lambda_{0}$ & 0 & 0 & 0 & 0\tabularnewline
\hline 
$\lambda_{1}$ & 0.821583899177118 & 0.821583899177077 & 0.821583899177230 & 0.821583899176988\tabularnewline
\hline 
$\lambda_{2}$ & 2.888537785769291 & 2.888537785769243  & 2.888537785769405 & 2.888537785769792\tabularnewline
\hline 
$\lambda_{3}$ & 2.944846615497959 & 2.944846615497851 & 2.944846615498256 & 2.977377367029730\tabularnewline
\hline 
$\lambda_{4}$ & 3.341726289664183 & 3.341726289664230 & 3.341726289664046 & 3.341726289664970\tabularnewline
\hline 
$\lambda_{5}$ & 4.550747949109708 & 4.550747949109686 & 4.550747949110111 & 4.550747949110250\tabularnewline
\hline 
$\lambda_{6}$ & 5.036739639826136 & 5.036739639826031 & 5.036739639826076 & 5.036739639826476\tabularnewline
\hline 
$\lambda_{7}$ & 6.233053526961343 & 6.233053526961285 & 6.233053526961188 & 6.233053526962100\tabularnewline
\hline 
$\lambda_{8}$ & 6.325490988924451 & 6.325490988924394 & 6.325490988924508 & 6.325490988924206\tabularnewline
\hline 
$\lambda_{\text{9}}$ & 7.805807719443767 & 7.805807719443299 & 7.805807719443544 & 7.805807719442640\tabularnewline
\hline 
$\lambda_{10}$ & 7.908416105951900 & 7.908416105952249 & 7.908416105952258 & 7.908416105952520\tabularnewline
\hline 
$\lambda_{11}$ & 9.404227647278619 & 9.404227647275778 & 9.404227647275357 & 9.404227647274947\tabularnewline
\hline 
\end{tabular}
\par\end{centering}
\caption{The first 12 eigenvalues $\lambda_{k},k=0,\ldots,11$  for different
numbers of grid points $N=2^{n},n=4,\ldots,10,12.$\label{tab:cosaini oval}}
\end{table}

\subsection{Optimization Solvers}

We solve the nonlinear system of ODEs (\ref{eq: a_k ODE}) in Section
\ref{subsec:Optimization-Solvers} by using the forward Euler method
with the time step $h$ to obtain the solution at $t+h.$ We can then
repeat this procedure iteratively until it finds the optimal shape.
To prevent the spurious growth of the high-frequency modes generated
by round-off error, we use 25th-order Fourier filtering and also filter
out the coefficients which is below $10^{-14}$ as used in \cite{nie1995stability}
after each iteration. 

\begin{figure}[H]
\begin{centering}
\subfloat[]{\begin{centering}
\includegraphics[scale=0.6]{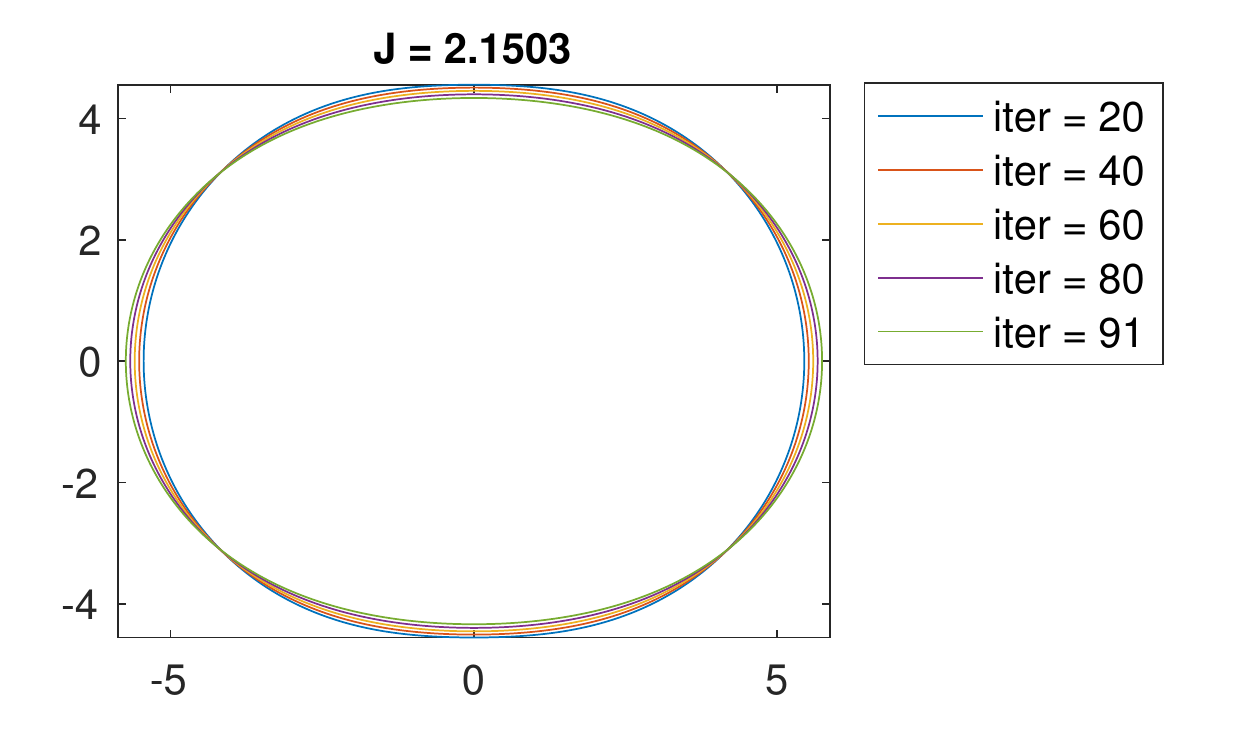}
\par\end{centering}
}\subfloat[]{\begin{centering}
\includegraphics[scale=0.6]{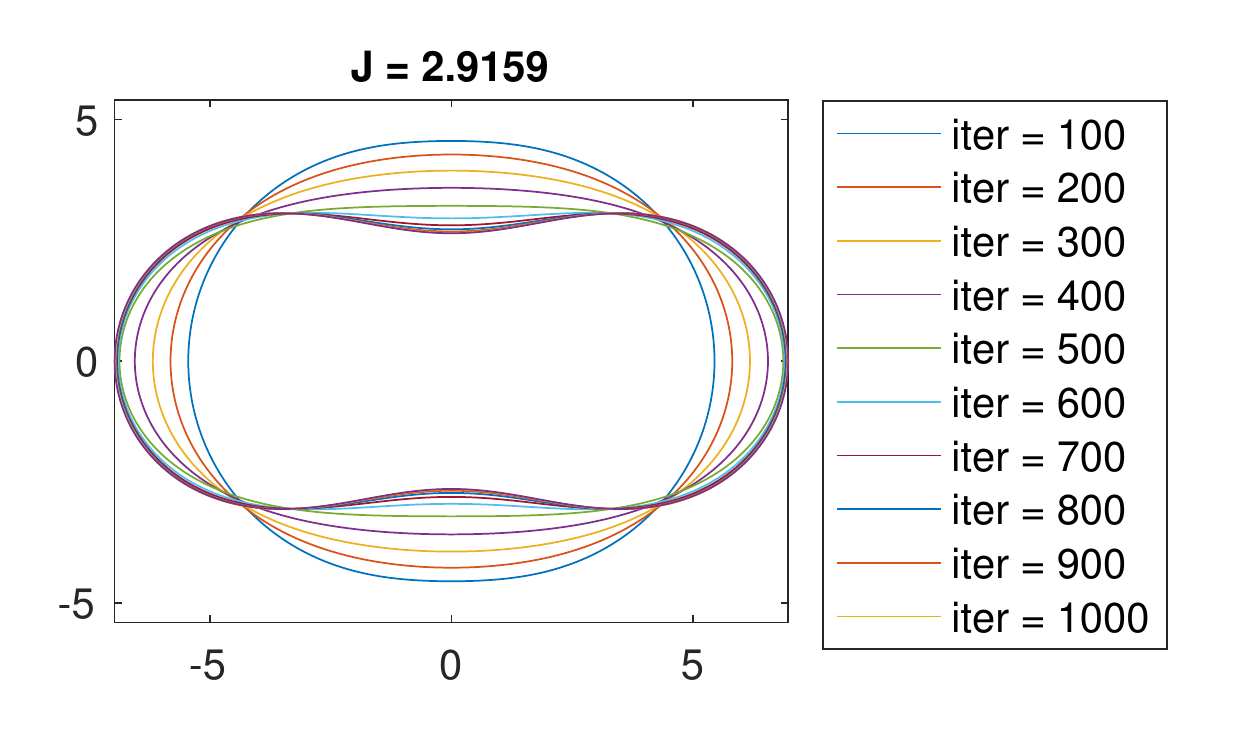}
\par\end{centering}
}
\par\end{centering}
\caption{Optimization of $\lambda_{2}^{A}(\Omega)$ without and with smoothing
$k$ are shown on (a) and (b), respectively. \label{fig:Optimization-of--1}}
\end{figure}

\begin{figure}[H]
\begin{centering}
\subfloat[]{\includegraphics[scale=0.15]{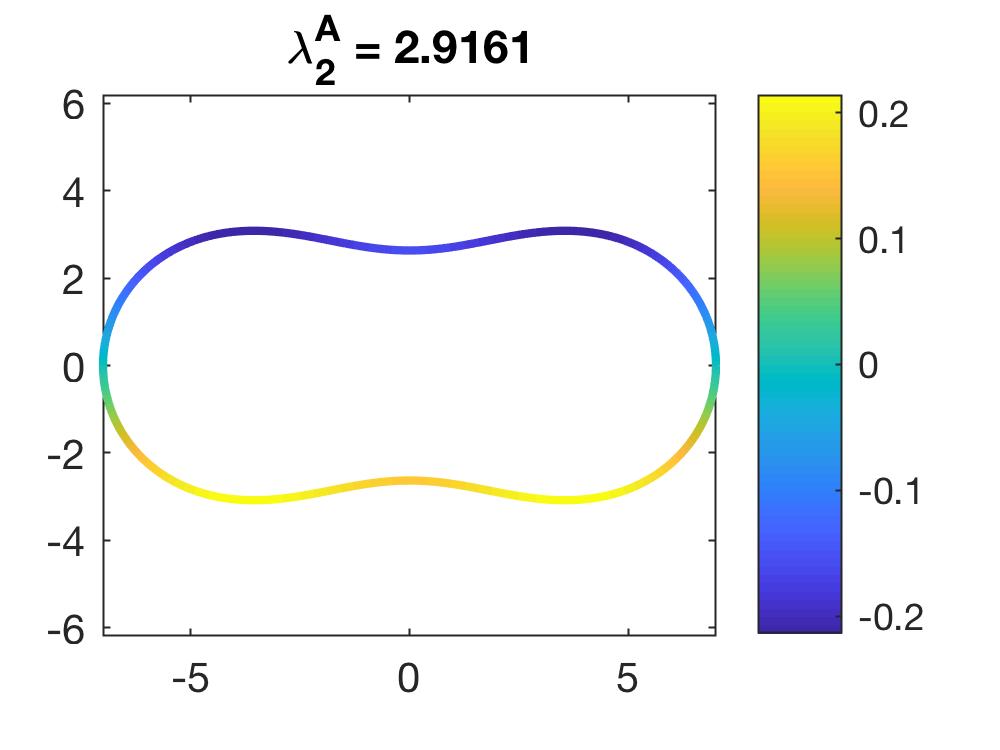}

}\subfloat[]{\includegraphics[scale=0.15]{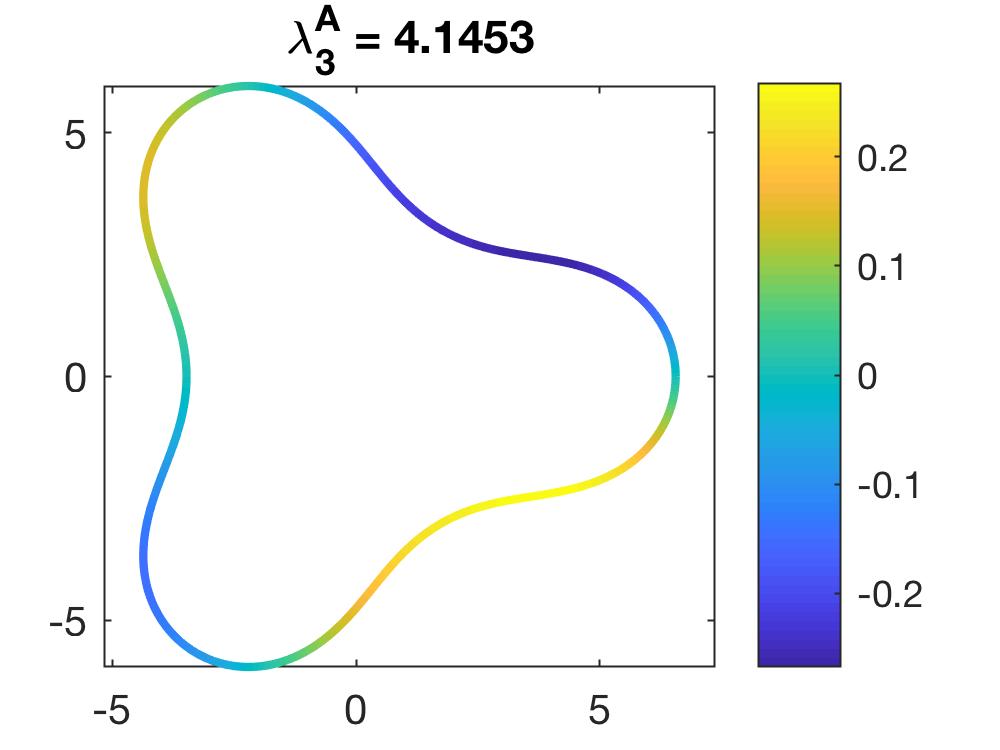}

}\subfloat[]{\includegraphics[scale=0.15]{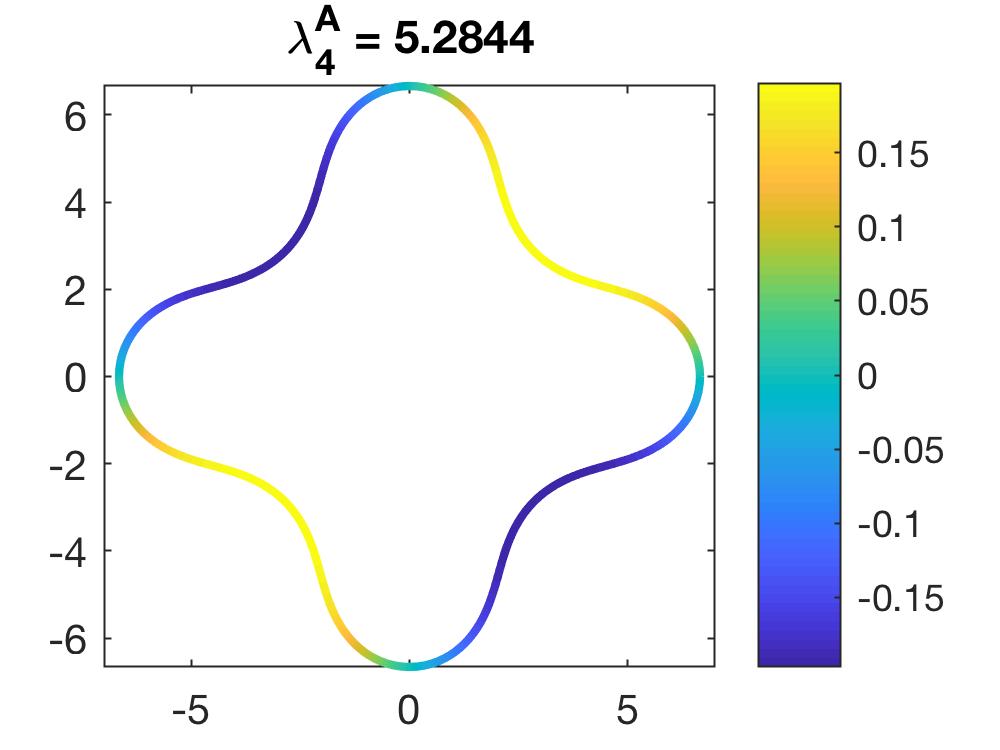}

}
\par\end{centering}
\begin{centering}
\subfloat[]{\begin{raggedright}
\includegraphics[scale=0.15]{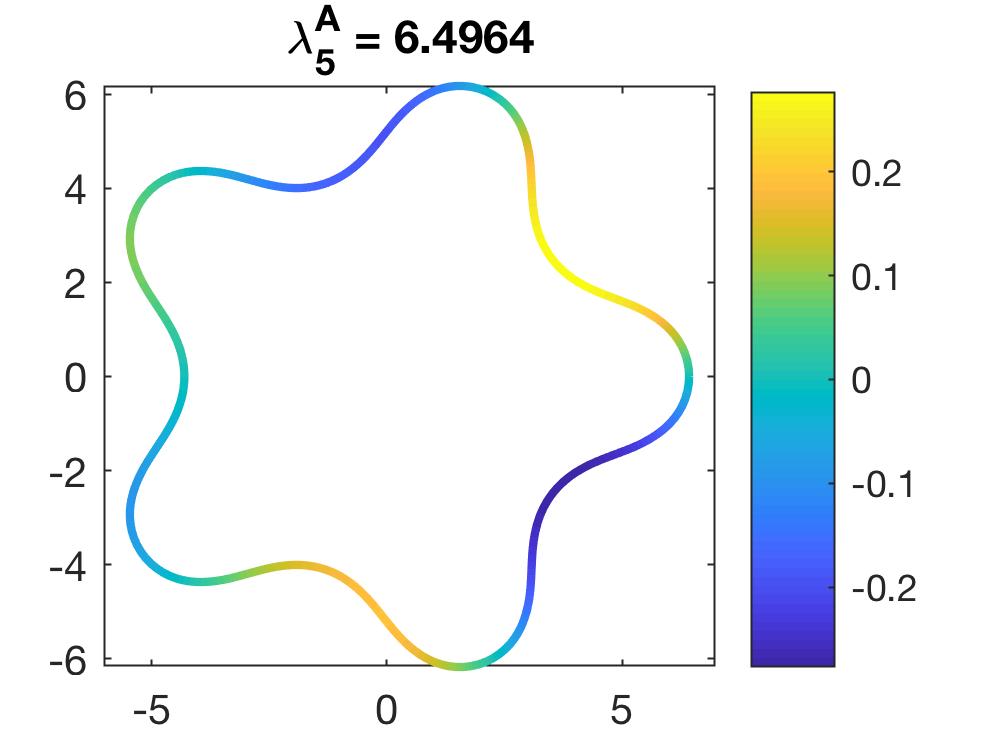}
\par\end{raggedright}

}\subfloat[]{\includegraphics[scale=0.15]{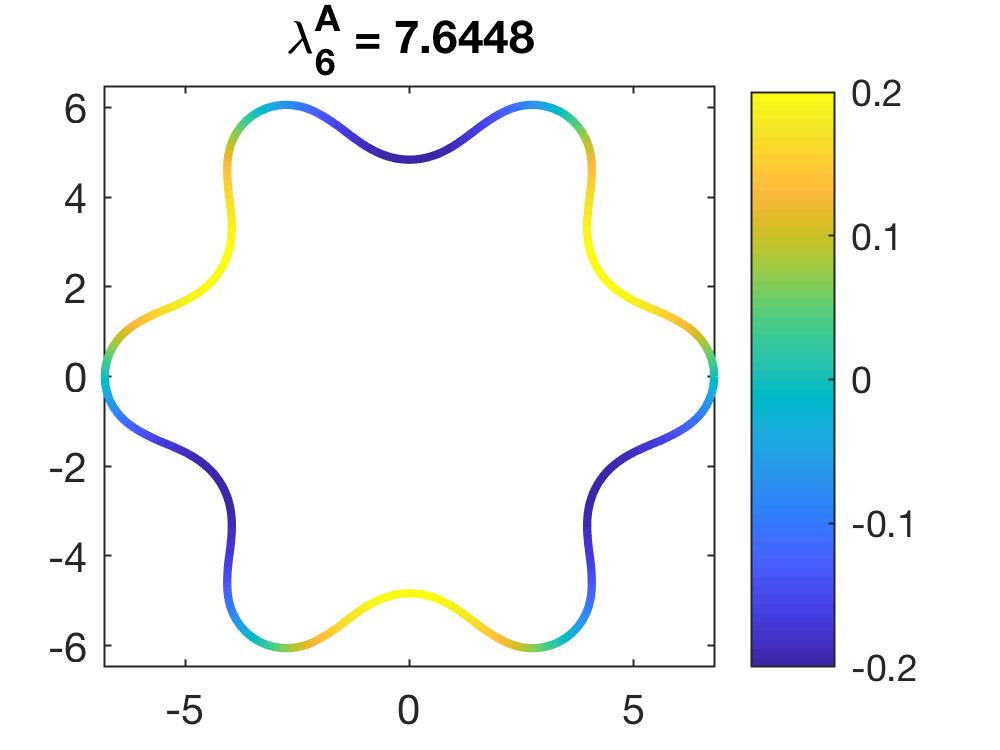}

}\subfloat[]{\begin{raggedright}
\includegraphics[scale=0.15]{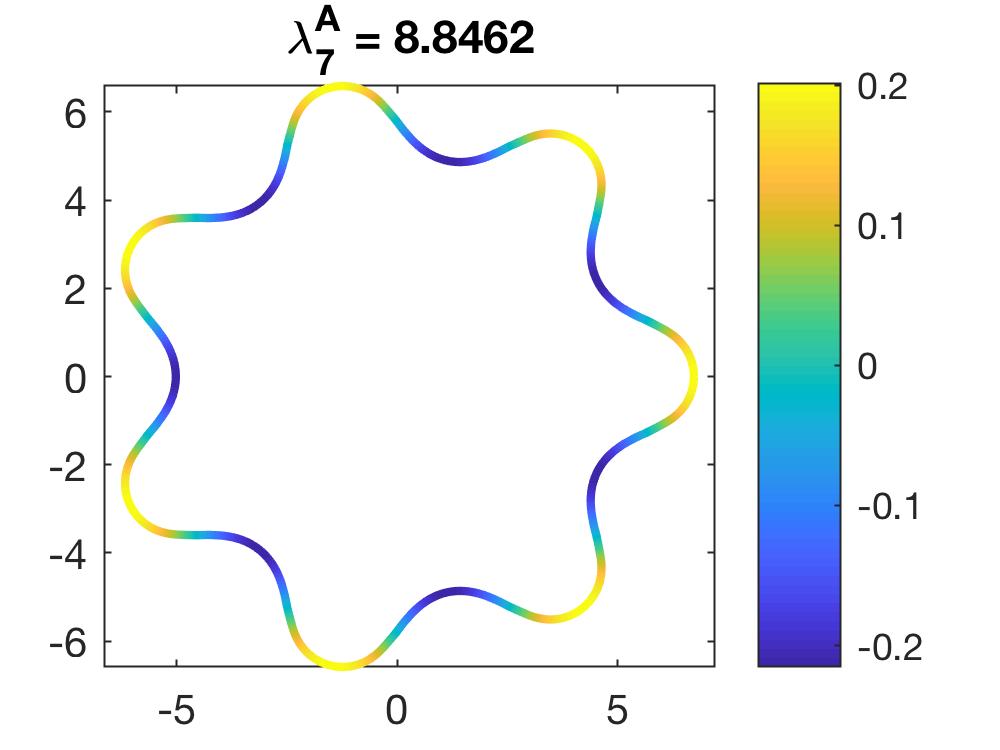}
\par\end{raggedright}

}
\par\end{centering}
\caption{The optimal shape of maximizing $\lambda_{n}^{A},n=2,\ldots,7.$\label{fig:The-optimal-shape}
The colors on the curve indicate the values of eigenfunctions. }

\end{figure}

In Figure \ref{fig:Optimization-of--1}(a), we show the evolution
of optimization of $\lambda_{2}^{A}$ with number of grid points $N=256$.
We start with a shape with a two-fold symmetry $f(w)=w+0.5w^{3}$
whose $\lambda_{2}^{A}=1.7791$. The algorithm was able to deform
the shape and increase the eigenvalue $\lambda_{2}^{A}$ up to 2.1503.
After that, the shape starts to generate kinks. Due to so-called crowding
phenomenon \cite{delillo1994accuracy}, the accuracy of the conformal
mapping will be effected and the shape will lose its smoothness. Thus,
we avoid this problem by smoothing the curvature term $\kappa$ in
the $z-$plane based on the moving average method with span $5$.
Using this smoothing technique at each iteration helps us to achieve
better results as shown in Figure \ref{fig:Optimization-of--1}(b).
In addition to smoothing, we also refine our time steps. We start
with an initial time step $h=0.1$ and halve the time step for every
time period $T=100$ and compute up to $5T.$ The optimal eigenvalues
$\lambda_{k}^{A},k=1,\cdots,7$ are summarized in Table \ref{tab:optimization-of-for}
and the optimal shapes which have $k$-fold symmetry are shown in
Figure \ref{fig:The-optimal-shape}. As observed in \cite{akhmetgaliyev2017computational},
the domain maximizing the $k-$th Steklov eigenvalue has $k$-fold
symmetry, and has at least one axis of symmetry. The $k-$th Steklov
eigenvalue has multiplicity 2 if $k$ is even and multiplicity 3 if
$k$ \ensuremath{\ge} 3 is odd. The first few nonzero coefficients
of the mapping function $f(w)$ of the optimal shapes are summarized
in Table \ref{tab:The-first-non} for $\lambda_{2}^{A}-\lambda_{7}^{A}$.
When optimizing $\lambda_{k}^{A}$, the optimal coefficients have
nonzero values for $a_{1+nk}$ where $n\in\mathbb{N}$. 

\section{Summary and Discussion}

We have developed a spectral method based on conformal mappings to
a unit circle to solve Steklov eigenvalue problem on general simply-connected
domains efficiently. Unlike techniques based on finite difference
methods or finite elements methods which requires discretization on
the general domains with boundary treatments, the method that we proposed
only requires discretization of the boundary of a unit circle. We
use a series expansion to represent eigenfunctions so that the discretization
leads to an eigenvalue problem for Fourier coefficients. In addition,
we study the maximization of area-normalized Steklov eigenvalue $\lambda_{k}^{A}$
based on shape derivatives and formulate this shape evolution in the
complex plane via the gradient ascent approach. With smoothing technique
and choices of time steps, we were able to find the optimal area-normalized
eigenvalues $\lambda_{k}^{A}$ for a given $k$. 

As aforementioned, the optimization of Steklov eigenvalue problems
on general non-simply-connected domains is a challenge open question.
This will require robust and efficient forward solvers of Steklov
eigenvalues and numerical techniques to perform shape optimizations
which may involve topological changes. In the near future, we plan
to explore the possibility in this direction by using Level Set approaches. 

\begin{table}[H]
\begin{centering}
\begin{tabular}{|c|c|c|}
\hline 
$\lambda_{2}^{A}$ & $\lambda_{3}^{A}$ & $\lambda_{4}^{A}$\tabularnewline
\hline 
\hline 
0 & 0 & 0\tabularnewline
\hline 
0.776986933500041 & 1.079861668314576 & 1.171320134341248\tabularnewline
\hline 
\textcolor{blue}{2.916071256633050} & 1.079861668314618 & 1.171320134341342\tabularnewline
\hline 
\textcolor{blue}{2.916071256753514} & \textcolor{blue}{4.145300664720734} & 1.611279604736676\tabularnewline
\hline 
3.277492771330297 & \textcolor{blue}{4.145300664720919} & \textcolor{blue}{5.284432268416950}\tabularnewline
\hline 
4.498623058633566 & \textcolor{blue}{4.145300672478222} & \textcolor{blue}{5.284433071016992}\tabularnewline
\hline 
5.041166283776032 & 4.914601402877488 & 5.448244774262810\tabularnewline
\hline 
6.118061463397883 & 6.024394262148678 & 5.448244774262829\tabularnewline
\hline 
6.272697585592614 & 6.024394262148718 & 6.489865254319582\tabularnewline
\hline 
7.693637484890079 & 7.628170417847103 & 7.335382999100261\tabularnewline
\hline 
7.809873534891437 & 7.628170417847109 & 7.335382999100267\tabularnewline
\hline 
9.262237946100434 & 8.953916828143468 & 8.636733197287754\tabularnewline
\hline 
\end{tabular}
\par\end{centering}
\begin{centering}
\begin{tabular}{|c|c|c|}
\hline 
$\lambda_{5}^{A}$ & $\lambda_{6}^{A}$ & $\lambda_{7}^{A}$\tabularnewline
\hline 
\hline 
0 & 0 & 0\tabularnewline
\hline 
1.239226322386241 & 1.265308570439713 & 1.291290525113730\tabularnewline
\hline 
1.239226322386290 & 1.265308570450229 & 1.291290525113829\tabularnewline
\hline 
1.945145428557867 & 2.117586845334797 & 2.250312782549877\tabularnewline
\hline 
1.945145428557917 & 2.117586845350534 & 2.250312782549968\tabularnewline
\hline 
\textcolor{blue}{6.496444238784153} & 2.427189796272854 & 2.777589136940805\tabularnewline
\hline 
\textcolor{blue}{6.496444238784278} & \textcolor{blue}{7.644759577423688} & 2.777589137507856\tabularnewline
\hline 
\textcolor{blue}{6.496444784959914} & \textcolor{blue}{7.644765127633966} & \textcolor{blue}{8.846228548846659}\tabularnewline
\hline 
6.732142619373287 & 7.771465908528654 & \textcolor{blue}{8.846229141938371}\tabularnewline
\hline 
6.732142619373318 & 7.771465908584982 & \textcolor{blue}{8.846229145378315}\tabularnewline
\hline 
8.128106267565293 & 7.979288943927369 & 9.050146643762274\tabularnewline
\hline 
8.803176067403113 & 7.979288943929372 & 9.050146643762337\tabularnewline
\hline 
\end{tabular}
\par\end{centering}
\caption{The optimization of $\lambda_{n}^{A},n=2,\ldots,7$ for the first
12 eigenvalues. \label{tab:optimization-of-for} }
\end{table}

\begin{table}[H]
\begin{centering}
\begin{tabular}{|c|c||c|c||c|c|}
\hline 
$\lambda_{2}^{A}$ &  & $\lambda_{3}^{A}$ &  & $\lambda_{4}^{A}$ & \tabularnewline
\hline 
\hline 
$a_{1}$ & 3.482625488377397 & $a_{1}$ & 4.172312832330094  & $a_{1}$ & 4.646184610628929\tabularnewline
\hline 
$a_{3}$ & 1.316760069380197 & $a_{4}$ & 1.018987204748702 & $a_{5}$ & 0.871201920168631\tabularnewline
\hline 
$a_{5}$ & 0.754288548863893 & $a_{7}$ & 0.514733681728398 & $a_{9}$ & 0.426363809028874\tabularnewline
\hline 
$a_{7}$ & 0.476336868618610 & $a_{12}$ & 0.301544250312563 & $a_{13}$ & 0.248913034789274\tabularnewline
\hline 
$a_{9}$ & 0.313178226238119 & $a_{15}$ & 0.187629804532372 & $a_{17}$ & 0.156013896896941\tabularnewline
\hline 
$a_{11}$ & 0.210225589908090 & $a_{18}$ & 0.120456778256507 & $a_{21}$ & 0.101460942653213\tabularnewline
\hline 
$a_{13}$ & 0.142829963279173 & $a_{21}$ & 0.078773331531670 & $a_{25}$ & 0.067443016703592\tabularnewline
\hline 
$a_{15}$ & 0.097776540001909 & $a_{24}$ & 0.052126631797892 & $a_{29}$ & 0.045465744838959\tabularnewline
\hline 
\end{tabular}
\par\end{centering}
\begin{centering}
\begin{tabular}{|c|c||c|c||c|c|}
\hline 
$\lambda_{5}^{A}$ &  & $\lambda_{6}^{A}$ &  & $\lambda_{7}^{A}$ & \tabularnewline
\hline 
\hline 
$a_{1}$ & 4.807404499929070 & $a_{1}$ & 5.298095057399003 & $a_{1}$ & 5.434176832482816\tabularnewline
\hline 
$a_{6}$ & 0.718033397997455 & $a_{7}$ & 0.665755972200186 & $a_{8}$ & 0.583992686042936\tabularnewline
\hline 
$a_{11}$ & 0.339254189743543 & $a_{13}$ & 0.310395425069731 & $a_{15}$ & 0.267954925737438\tabularnewline
\hline 
$a_{16}$ & 0.195019993266578 & $a_{19}$ & 0.178491217642749 & $a_{22}$ & 0.153351975823925\tabularnewline
\hline 
$a_{21}$ & 0.121279950688959 & $a_{25}$ & 0.111622289351001 & $a_{29}$ & 0.095872749417195\tabularnewline
\hline 
$a_{26}$ & 0.078576438618779 & $a_{31}$ & 0.072927620185693 & $a_{36}$ & 0.062767599076665\tabularnewline
\hline 
$a_{31}$ & 0.052167793728917 & $a_{37}$ & 0.048905352581595 & $a_{43}$ & 0.042229353541252\tabularnewline
\hline 
$a_{36}$ & 0.035185234307276 & $a_{43}$ & 0.033346733965258 & $a_{50}$ & 0.028892884585296\tabularnewline
\hline 
\end{tabular}
\par\end{centering}
\caption{The first few nonzero coefficients of the mapping function $f(w)$
of the optimal shapes for $\lambda_{2}^{A}-\lambda_{7}^{A}$.\label{tab:The-first-non}}
\end{table}


\bibliography{SteklovOptim}

\end{document}